\documentclass[3p,times]{elsarticle}
\usepackage{amsmath,amssymb,amsfonts}
\usepackage{amsthm}
\usepackage{mathrsfs}

\usepackage{graphicx}
\usepackage{subcaption}

\usepackage{booktabs}
\usepackage{multirow}

\usepackage[linesnumbered,ruled,vlined]{algorithm2e}
\usepackage{xcolor}

\usepackage{hyperref}

\usepackage{booktabs}

\usepackage{siunitx}

\newtheorem{theorem}{Theorem}[section]

\theoremstyle{definition}

\theoremstyle{remark}
\newtheorem{remark}[theorem]{Remark}

\newcommand{\bx}{\mathbf{x}}

\journal{Journal of Computational Physics}

\usepackage{listings}%

\usepackage{graphicx} 			%图片
\usepackage{subcaption}			%子图
\usepackage{float}

\newcommand{\inpd}[2]{\left\langle #1, #2 \right\rangle}

\usepackage{tikz}
\usetikzlibrary{positioning,arrows.meta,shapes.geometric,calc}
\usepackage{amsmath}

\usepackage{verbatim}

\usepackage{lineno}

\begin{document}
%\linenumbers
\begin{frontmatter}

\title{A Deep Ritz Method for High-Dimensional Steady States of the Cahn--Hilliard Equation}

\author[1]{Yi Liu }
\author[1]{Shuting Gu\textsuperscript{✳} }
\cortext[cor1]{Corresponding author. Email: gushuting@sztu.edu.cn} 
\address[1]{School of Artificial Intelligence, Shenzhen Technology University, Shenzhen 518118, P.R. China }
%\address[2]{School of Artificial Intelligence, Shenzhen Technology University, Shenzhen 518118, P.R. China }

\begin{abstract}
%%%
The Cahn--Hilliard equation is a fundamental model for describing phase separation phenomena in binary mixtures. Traditional numerical methods, such as finite difference and finite element methods, often incur substantial computational cost, particularly when computing steady-state solutions in high-dimensional settings.
To address this challenge, we propose a deep learning-based framework-- the Deep Ritz method--for computing steady states of the Cahn--Hilliard equation under periodic boundary conditions. An enhanced augmented Lagrangian formulation is incorporated to strictly enforce the mass conservation constraint, while separable Fourier feature mappings are employed to naturally encode periodicity and enhance the representation of nontrivial solution structures.
The proposed method exhibits a notable dual capability: it not only achieves fast convergence to steady states but also effectively identifies multiple nontrivial solutions corresponding to different local minimizers of the energy functional. Extensive numerical experiments in one-, two-, and three-dimensional cases demonstrate that the method can successfully capture a rich variety of phase separation patterns, including droplet-type, lamellar, and tubular structures. We also compare the numerical results with those computed by finite difference method (FDM) for one- and two-dimensional cases, highlighting the effectiveness and robustness of the proposed approach in exploring complex high-dimensional energy landscapes.

%%%%
\end{abstract}


\end{frontmatter}

%\linenumbers

%% main text

\section{Introduction}

The Cahn--Hilliard (CH) equation is a fundamental partial differential equation for modeling phase separation processes in binary mixtures, polymer blends, and multiphase systems \citep{cahn1958free}. It describes the evolution of a conserved order parameter driven by the minimization of a free energy functional. Computing steady-state solutions of the CH equation is of particular importance, as they correspond to equilibrium configurations of the underlying physical system.
Traditional numerical methods, such as the finite difference and finite element methods, have been widely used to solve the CH equation. However, these approaches typically rely on mesh-based discretizations, and their computational cost grows rapidly with the spatial dimension $d$, suffering from the well-known curse of dimensionality \cite{elliott1989global}. This limitation makes high-dimensional steady-state computations particularly challenging.

In recent years, deep learning-based methods have emerged as a promising alternative to solve high-dimensional partial differential equations \cite{han2017deep, han2018solving, raissi2019physics, beck2019machine}. These methods leverage the expressive power of neural networks for high-dimensional function approximation and have demonstrated significant potential in mitigating the curse of dimensionality. 
This is supported by the universal approximation theorem \cite{hornik1991approximation}. Moreover, it has been shown in \cite{zhang2022deep} that, with appropriate choices of activation functions, a fully connected feedforward neural network with width $36d(2d+1)$ and depth $11$ is capable of approximating any continuous function on a hypercube of dimensions $d$ to arbitrary accuracy.

 In general, deep learning-based methods can be broadly classified into two categories: direct approaches derived from PDE formulations and stochastic approaches originating from probabilistic models.

 %% the follwoing is to be simplified.

In the first category, deep neural networks (DNNs) are employed to solve PDEs directly. For instance, physics-informed neural networks (PINNs) were introduced in \cite{raissi2017physics, raissi2019physics}, where the loss function incorporates observational data as well as initial and boundary conditions to approximate PDE solutions. The deep Galerkin method (DGM) was proposed in \cite{sirignano2018dgm} to solve high-dimensional free-boundary parabolic PDEs. In addition, the Deep Ritz method \cite{weinan2018deep} was developed for variational problems arising from PDEs. The weak adversarial network (WAN) approach \cite{zang2020weak} leverages weak formulations to handle general high-dimensional PDEs defined on complex domains.
Furthermore, a deep feedforward network, called PDE-Net \cite{pmlr-v80-long18a}, was proposed to learn and predict the dynamics of complex systems. Its extension, PDE-Net 2.0 \cite{long2019pde}, further enhances flexibility and expressive power by simultaneously learning differential operators and nonlinear response functions of the underlying PDE models.

In the second category, DNNs are employed to solve PDEs via probabilistic representations. The deep BSDE method was proposed for a class of nonlinear PDEs based on forward–backward stochastic differential equations (FBSDEs) \cite{han2017deep, han2018solving, han2025brief}. Its variants have been further developed in \cite{raissi2024forward, zhang2022fbsde, ji2020three}, extending this framework through pathwise comparisons of stochastic processes. This methodology has also been applied to stochastic control problems \cite{han2016deep}.
Building upon the BSDE representation of PDEs, the backward deep learning (DBDP) algorithm \cite{hure2020deep} and the multistep deep backward dynamic programming (MDBDP) scheme \cite{germain2022approximation} were proposed to analyze approximation errors within the framework of neural network approximation theory. The MDBDP, in particular, provides an error analysis for the deep splitting scheme applied to semilinear PDEs \cite{beck2021deep}.
In addition to BSDE-based approaches, another promising direction is the martingale formulation inspired by Varadhan’s martingale problem. This perspective offers a new avenue for tackling high-dimensional PDEs and stochastic control problems \cite{cai2023deepmartnet}. In this context, the martingale-based deep neural network (DeepMartNet) \cite{cai2023deepmartnet} has been developed to approximate solutions to boundary value problems and eigenvalue problems for elliptic PDEs.

% up to here. the uper is to be simplified.

%Representative approaches include physics-informed neural networks (PINNs) \cite{raissi2019physics}, the deep Galerkin method (DGM) \cite{sirignano2018dgm}, and the Deep Ritz method \cite{weinan2018deep}, among others.

Despite their success, several challenges remain when applying deep learning methods to phase-field models such as the CH equation. In particular, the accurate enforcement of physical constraints, such as mass conservation, and the treatment of periodic boundary conditions are nontrivial. Moreover, identifying multiple nontrivial steady states corresponding to different local minimizers of the energy functional remains largely unexplored.

In this work, we develop a deep learning-based framework for computing steady states of the CH equation in high-dimensional settings. The proposed method is built upon the Deep Ritz formulation \cite{weinan2018deep} and incorporates an augmented Lagrangian approach to strictly enforce mass conservation. To handle periodic boundary conditions, we introduce separable Fourier feature mappings \cite{tancik2020fourier, rahaman2019spectral}, which naturally encode periodicity and enhance the representation of high-frequency solution structures.

The main contributions of this work include the following aspects: (1) a Deep Ritz-based method is proposed for computing steady states of the CH equation in high-dimensional domains, alleviating the curse of dimensionality. (2) an augmented Lagrangian formulation is incorporated to effectively enforce mass conservation within the neural network framework. (3) the Fourier feature mapping technique is introduced to naturally handle periodic boundary conditions and improve the approximation of complex solution patterns.
 This proposed method is capable of identifying multiple nontrivial steady states associated with different local minimizers of the Ginzburg-Landau free energy functional.
 We test extensive numerical experiments in one-, two-, and three-dimensional settings, which demonstrate the accuracy, efficiency, and robustness of the method.

The remainder of the paper is organized as follows. In Section~\ref{Review}, we briefly review the Deep Ritz method. Section~\ref{Main_method} presents the augmented Lagrangian formulation, the neural network architecture, the Fourier feature mapping technique and the associated algorithm. Numerical results on one-, two-, and three-dimensional CH equations are reported in Section~\ref{Numerical_results}, followed by concluding remarks, extensions and future research directions in Section~\ref{Conclusion}.

\section{Review of the Deep Ritz Method}\label{Review}

The Deep Ritz Method (DRM), introduced in \cite{weinan2018deep}, is a neural-network-based approach for solving variational problems or elliptic PDEs. In DRM, trial functions are parameterized by deep neural networks, and the associated energy functional is minimized directly. To illustrate, consider the classical Poisson problem:
\begin{align}\label{Ellip_PDE}
    \begin{cases}
        - \Delta u = f, \quad\text{in } \Omega, \\
        \quad u=0, \quad \text{on } \partial\Omega.
    \end{cases}
\end{align}
The solution of \eqref{Ellip_PDE} can be equivalently formulated as the minimizer of the variational problem
$$ \min\limits_{u\in H} E(u), $$
with
\[
\displaystyle E(u) = \int_{\Omega} \displaystyle \frac{1}{2} \big|\nabla u(\bx)\big|^2 - f(\bx) \, u(\bx) d\mathbf{x},
\]
where $H$ denotes the space of admissible trial functions.

In DRM, a deep neural network \( u_\theta(\mathbf{x}):= u(\bf x;\theta) \)  with parameters \(\theta\) is employed to approximate the unknown solution \(u(\mathbf{x})\). Boundary conditions are typically enforced via a penalty term, yielding the following loss function:
\begin{equation}
\mathcal{L}(\theta) = \int_{\Omega} \left(\frac{1}{2} \big| \nabla u_\theta({\bf x}; \theta) \big|^2 - f({\bf x}) \, u_\theta({\bf x}; \theta) \right)d\mathbf{x}
+ \nu\,\int_{\partial\Omega} \big|u_\theta({\bf x}; \theta) \big|^2\,dS    
\end{equation}
where $\nu$ is a sufficiently large penalty parameter. The original PDE problem is then reduc0.0ed to the following optimization problem:
\begin{equation}
    \min\limits_{\theta} \mathcal{L}(\theta),
\end{equation}
which is typically solved numerically using mini-batch stochastic gradient descent (SGD).

The key advantages of DRM include its mesh-free formulation, compatibility with stochastic mini-batch training, and potential scalability to moderately high dimensions. Subsequent theoretical analysis have provided generalization bounds for DRM, demonstrating that, under appropriate assumptions, the approximation error does not necessarily grow with the spatial dimension \(d\) \cite{Lu2021Generalization}.  

\section{Problem Setup and Main Method}\label{Main_method}

In this work, we focus on the steady-state solutions of the Cahn–Hilliard equation \cite{chaikin2000principles} under periodic boundary conditions, 
\begin{equation}\label{CH}
    \partial_t \phi = \Delta \delta_\phi F(\phi) = -\varepsilon^2 \Delta^2 \phi + \Delta (\phi^3 -\phi),
\end{equation}
which corresponds to the $H^{-1}$-gradient flow of the Ginzburg–Landau free-energy functional,
\[
F(\phi) = \int_{\Omega} \left( \frac{\varepsilon^2}{2} |\nabla \phi|^2 + \frac14(\phi^2 - 1)^2 \right)d\mathbf{x},
\]
where $\varepsilon>0$ is a mobility parameter and $\phi(\mathbf{x})$ is an order parameter representing the concentration of one component in a binary alloy.

Most existing approaches locate the steady states of the CH equation by directly solving the time-dependent PDE \eqref{CH} using traditional numerical methods. However, such approaches are often computationally prohibitive due to the curse of dimensionality. In this work, we instead consider the problem from a variational perspective, seeking minimizers of the Ginzburg–Landau free energy in the $H^{-1}$ metric. Since mass is conserved in the $H^{-1}$ framework, we employ an enhanced augmented Lagrangian method to explicitly enforce this constraint while searching for steady-state solutions.

\subsection{Enhanced Augmented Lagrangian Method}

In this work, we solve the constrained variational problem using a deep learning approach. Specifically, the order parameter $\phi(\mathbf{x})$
 is parameterized by a deep neural network $u_\theta(\mathbf{x})$
 with parameters $\theta$, i.e.,
 $$\phi(\mathbf{x}) \approx u_\theta(\mathbf{x}).$$

 To enforce the mass conservation constraint within the deep learning framework, we adopt an enhanced augmented Lagrangian (AL) method. This approach employs a nested (outer-inner) loop structure to simultaneously optimize the free energy and enforce the mass constraint. The corresponding loss functional is defined as
 \begin{equation}\label{Loss_AL}
   \mathcal{L}_{\mathrm{AL}}(\theta) 
= \underbrace{\int_{\Omega} \left( \frac{\varepsilon^2}{2} |\nabla u_\theta|^2 + \frac14(u_\theta^2 - 1)^2 \right)d\mathbf{x}}_{\text{energy term}} \quad 
+ \underbrace{\lambda ~ \Bigl( \bar u_\theta - m_0\Bigr)
+ \frac{\mu}{2} ~ \Bigl( \bar u_\theta - m_0\Bigr)^2 }_{\text{augmented Lagrange term}} := \mathcal{L}_{\mathrm{i}} + \mathcal{L}_{\mathrm{m}}, 
\end{equation}
where \(\lambda\) and \(\mu\) are positive parameters associated with the Lagrange multiplier and penalty term, respectively, and $m_0$ denotes the initial mass, i.e.,
$$\int_\Omega \phi({\bf{x}}) d \mathbf{x} = m_0.$$ 
 $\bar u_\theta$ is estimated via stochastic sampling (e.g., Sobol quasi-Monte Carlo points): 
\begin{equation}\label{ubar}
    \bar u_\theta = \int_{\Omega} u_\theta ~ d \mathbf{x} \approx \frac{1}{N}\sum_{i=1}^N u_\theta(\mathbf{x}_i).    
\end{equation}

A natural approach to enforce periodic boundary conditions is to incorporate them as an additional term in the loss function. For example, in one dimension with $\Omega = [a,b]$, the boundary loss is defined as
\begin{equation}
     \mathcal{L}_b(\theta) = \big|u_\theta(a) - u_\theta(b) \big|^2.
\end{equation}
The total loss is then given by 
$$ \mathcal{L}_{\text{total}}(\theta) = \mathcal{L}_{\mathrm{AL}}(\theta) + \alpha \mathcal{L}_b(\theta) = \mathcal{L}_{\mathrm{i}} + \mathcal{L}_{\mathrm{m}} + \alpha \mathcal{L}_{\mathrm{b}}, $$
where $\alpha>0$ is a weighting coefficient. The neural network parameters $\theta$ are updated by minimizing 
$\mathcal{L}_{\text{total}}(\theta)$, and the resulting minimizer corresponds to the steady-state solution of the CH equation.

\begin{remark}
The simple penalty-based approach for enforcing periodic boundary conditions is generally only capable of locating the trivial solution. For non-trivial solutions, this method often fails due to the spectral bias of standard neural networks. To overcome this limitation, we employ the Fourier Feature Mapping technique to handle periodic boundary conditions, as detailed in Subsection \ref{sec:fourier-mapping}.
\end{remark}

\subsection{Global Residual Deep Neural Network}

To solve the variational minimization problem within the deep learning framework, we adopt a global residual network architecture, referred to as the \textit{Global Residual Deep Neural Network} (GlobalResNet), to enhance both stability and accuracy. The order parameter \(\phi(\mathbf{x})\) is parameterized by a neural network \(u_\theta(\mathbf{x})\).  

\begin{figure}[H]
    \centering
    \includegraphics[width=0.35\linewidth]{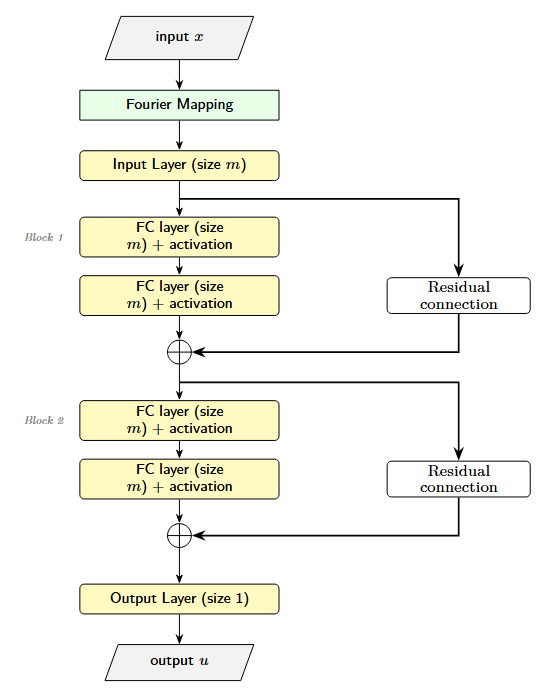}
   \caption{Schematic illustration of the GlobalResNet architecture. The model combines a global linear mapping with a nonlinear residual network composed of stacked residual blocks, enabling efficient learning of nontrivial solution structures.}
    \label{fig:global-res-deep-ritz}
\end{figure}

In this architecture, the network output is expressed as the sum of a global linear mapping and a nonlinear transformation learned through residual blocks (as illustrated in Figure \ref{fig:global-res-deep-ritz}). Specifically, for an input \(\mathbf{x} \in \mathbb{R}^d\), the network output is defined as
\begin{equation}
    u_\theta(\mathbf{x}) = \beta_1 \, \mathbf{w}^\top \mathbf{x} + \beta_2 + \mathcal{N}_\theta(\mathbf{x}),
\end{equation}
where \(\beta_1, \beta_2 \in \mathbb{R}\) are trainable scalar parameters, \(\mathbf{w}\) is a learnable weight vector, and \(\mathcal{N}_\theta(\mathbf{x})\) represents the nonlinear component constructed by stacked residual blocks. In this work, the activation function is chosen as a ReLU-type nonlinearity, namely \(\max\{\mathbf{x}^{1.5}, 0\}\).  

The global residual structure provides two key advantages. First, it mitigates the tendency of the network to converge to trivial or nearly constant solutions, a common issue in physics-informed neural network training. Second, it accelerates convergence by learning the target function as a perturbation of a simple baseline (the linear mapping), rather than from scratch. Numerical experiments indicate that the proposed GlobalResNet achieves improved accuracy and generalization compared to standard feedforward networks or local residual architectures.

\begin{remark}
A common challenge in training deep neural networks for variational problems is the tendency to converge to trivial or nearly constant solutions, particularly when the optimization landscape contains flat regions. One key factor in avoiding such behavior is the choice of weight initialization.  

In this work, we adopt the Kaiming initialization scheme \cite{He_2015_ICCV}, where the weights are initialized according to
\[
W_{ij} \sim \mathcal{N}\!\left(0, \frac{2}{n_{\text{in}}}\right),
\]
with \(n_{\text{in}}\) denoting the number of input neurons. This initialization preserves the variance of activations across layers when using ReLU-type nonlinearities, thereby mitigating vanishing or exploding signal issues.  

From a practical perspective, Kaiming initialization enhances the expressive capacity of the network during the early stages of training, allowing it to explore richer solution landscapes and avoid collapse to trivial constant solutions. As a result, it facilitates the discovery of nontrivial steady-state solutions while maintaining training stability and convergence efficiency.
\end{remark}

\subsection{Fourier Feature Mapping}
\label{sec:fourier-mapping}
Coordinate-based neural networks (also known as coordinate MLPs) have been widely used to approximate solutions of partial differential equations. However, such models exhibit a well-known limitation, namely the \emph{spectral bias}, whereby low-frequency components are learned preferentially while high-frequency or sharply varying features (e.g., narrow interfacial layers) are difficult to capture \cite{Tancik2020FourierFeatures}.  

To alleviate this issue, we employ the Fourier Feature Mapping (FFM) technique, which lifts the input coordinates into a higher-dimensional feature space before feeding them into the neural network. Specifically, the mapping $\gamma(\mathbf{x})$ is defined as
\begin{equation}
\gamma(\mathbf{x}) : \mathbf{x} \in \mathbb{R}^d \mapsto 
\begin{bmatrix}
\cos \bigl(2\pi\,\mathbf{B}\mathbf{x}\bigr)\\[2pt]
\sin \bigl(2\pi\,\mathbf{B}\mathbf{x}\bigr)
\end{bmatrix}
\in \mathbb{R}^{2m},
\end{equation}
where $\mathbf{B} \in \mathbb{R}^{m \times d}$ is a frequency matrix whose rows correspond to frequency vectors. This transformation enriches the function space accessible to the network and significantly improves its ability to represent high-frequency components. Moreover, by appropriately choosing integer-valued frequencies, the mapping can enforce exact periodicity along specified coordinate directions.

A commonly used variant is the Random Fourier Feature (RFF) mapping, where the entries of $\mathbf{B}$ are sampled from a Gaussian distribution, i.e., $b_{ij} \sim \mathcal{N}(0, s^2)$. The parameters $m$ and $s$ control the number and magnitude of the frequencies, respectively, thereby determining the representational capacity of the mapping. While RFF provides a flexible and effective way to mitigate spectral bias \cite{Tancik2020FourierFeatures, Rahimi2007RandomFeatures}, the sampled frequencies are generally non-integer and thus do not guarantee exact periodicity at the domain boundaries.

For problems defined on periodic domains, it is advantageous to employ integer-frequency Fourier bases. For example, the classical Fourier series expansion in one dimension,
\[
\phi(x) \approx a_0 + \sum_{n=1}^{f_m} \Bigl[a_n \cos \bigl(2\pi n x / L_x\bigr) + b_n \sin \bigl(2\pi n x / L_x\bigr)\Bigr],
\]
can be viewed as a special case of Fourier feature mapping with integer frequencies. Here, $L_x$ denotes the period along the $x$-direction and $f_m$ is the maximum frequency.

For a $d$-dimensional periodic domain $\Omega = \prod_{j=1}^d [0, L_j]$, a natural extension is the \emph{Cartesian Fourier feature} construction. Given a maximum frequency $f_m$, the full frequency set is defined as the tensor product
\[
\mathcal{K}_f = \{-f_m, \dots, f_m\}^d \subset \mathbb{Z}^d.
\]
The frequency matrix $\mathbf{B} \in \mathbb{Z}^{m \times d}$ is formed by stacking all frequency vectors $\mathbf{n}_k = \{n_{kj}\} \in \mathcal{K}_f $, i.e., 
$\mathbf{B} = [\mathbf{n}_1; \mathbf{n}_2; \cdots; \mathbf{n}_m] $, 
resulting in
\[
m = (2f_m+1)^d.
\]

For non-unit domains, each frequency component is normalized by the corresponding period, i.e., $n_{kj} \mapsto n_{kj} / L_j, ~ j=1,\cdots, d, k = 1,\cdots, m$. This construction ensures exact periodicity along each coordinate direction.

However, the Cartesian construction suffers from an exponential growth in the number of modes with respect to the dimension $d$, leading to significant computational and memory costs. In addition, many high-frequency mixed modes may contribute little to the representation of the solution, especially for problems with moderate smoothness.

To address this issue, we adopt a more efficient alternative, namely the \emph{separable Fourier feature} construction. Instead of using the full tensor product, only axis-aligned frequency modes are retained. The corresponding frequency set is defined as
\[
\mathcal{K}_f^{\mathrm{sep}} 
= \{\mathbf{0}\} \cup \bigcup_{j=1}^d \left\{ (0,\dots,0,n_j,0,\dots,0) : n_j \in \{-f_m,\dots,f_m\}\setminus\{0\} \right\}.
\]
The total number of modes in this case is
\[
m_{\mathrm{sep}} = 1 + d \cdot (2f_m),
\]
which scales linearly with the dimension $d$. The corresponding frequency matrix $\mathbf{B}$ is constructed from vectors in $\mathcal{K}_f^{\mathrm{sep}}$.

The separable Fourier feature mapping preserves exact periodicity while significantly reducing the number of features. As a result, it alleviates the curse of dimensionality inherent in the full Cartesian construction, providing an efficient and scalable representation for high-dimensional problems.

\begin{remark}
The separable Fourier feature construction omits mixed-direction modes involving multiple nonzero frequency components (e.g., $(n_1,n_2)$ with $n_1 \neq 0$ and $n_2 \neq 0$), and therefore cannot directly represent such cross-frequency interactions. However, these mixed modes can be implicitly recovered through nonlinear interactions within the network.  

For instance, consider
\[
a(x) = \cos(2\pi n_1 x), \qquad b(y) = \cos(2\pi n_2 y).
\]
By the trigonometric identity,
\[
a(x)\, b(y)
= \cos(2\pi n_1 x)\cos(2\pi n_2 y)
= \frac{1}{2} \Bigl[ \cos\bigl(2\pi(n_1 x + n_2 y)\bigr) + \cos\bigl(2\pi(n_1 x - n_2 y)\bigr) \Bigr],
\]
which contains mixed-frequency components corresponding to the Cartesian modes $(n_1,n_2)$ and $(n_1,-n_2)$.  

Motivated by this observation, for periodic PDEs with fine-scale or high-frequency structures, we adopt a hybrid feature representation that combines:  
(i) separable integer-frequency sine/cosine features along each coordinate direction, which ensure exact periodic boundary conditions; and  
(ii) Random Fourier Features (RFF) with prescribed size $m$ and scale $s$, which enhance high-frequency expressivity.  This integrated construction preserves exact periodicity while providing richer spectral coverage and improved approximation capability. This design strikes a balance between computational efficiency and expressive power in high-dimensional settings.
\end{remark}

\begin{remark}
Fourier feature mappings introduce frequency-dependent scaling factors in the derivatives, which may amplify gradient magnitudes. However, this effect does not impact the optimization direction, and the feature values remain bounded. In practice, the method remains stable and converges to physically meaningful solutions, while benefiting from improved representation of high-frequency structures.
\end{remark}

\subsection{Algorithmic Framework}

With the incorporation of the Fourier feature mapping, periodic boundary conditions are satisfied intrinsically through the input representation. As a result, no additional boundary penalty term is required, and the total loss function reduces to the augmented Lagrangian functional defined in \eqref{Loss_AL}.

Specifically, the total loss is given by
\begin{equation}\label{L_tot_FFM}
\mathcal{L}_{\mathrm{total}}(\theta) 
= \int_{\Omega} \left( \frac{\varepsilon^2}{2} |\nabla u_\theta|^2 + \frac{1}{4}(u_\theta^2 - 1)^2 \right)\, d\mathbf{x}
+ \lambda \left( \bar u_\theta - m_0 \right)
+ \frac{\mu}{2} \left( \bar u_\theta - m_0 \right)^2,
\end{equation}
where $m_0$ denotes the prescribed mass, and $\lambda$ and $\mu$ are the Lagrange multiplier and penalty parameter, respectively, as defined in \eqref{Loss_AL}.

To solve the constrained minimization problem, we employ a nested optimization framework consisting of inner iterations and outer cycles, as summarized in Algorithm \ref{Algorithm}. 

In the inner iteration, the network parameters $\theta$ are updated by minimizing the total loss $\mathcal{L}_{\mathrm{total}}$, while keeping the Lagrange multiplier $\lambda$ and the penalty parameter $\mu$ fixed.

In each outer cycle, the Lagrange multiplier $\lambda$ and the penalty coefficient $\mu$ are updated
according to:
\[
    \lambda^{(k+1)} = \lambda^{(k)} + \mu^{(k)} \bigl(\bar{u}_\theta - m_0\bigr), \quad
    \mu^{(k+1)} = \min\bigl(\rho\, \mu^{(k)}, \mu_{\max}\bigr),
\]
where $\rho > 1$ denotes the growth rate of the penalty parameter, and $\mu_{\max}$ is a prescribed upper bound. Here, $\bar{u}_\theta$, the mean value of the network output, is used to update the Lagrange multiplier. 

To enhance stability and convergence in practice, the Lagrange multiplier $\lambda$ is initialized as zero in the early outer iterations, allowing the penalty term $\mu$ to dominate the enforcement of the mass constraint;  As the training proceeds, $\lambda$ gradually becomes active and refine the mean value $\bar{u}_\theta$ toward the prescribed mass $m_0$.
This nested optimization strategy provides an efficient and stable mechanism for enforcing constraints while maintaining fast convergence in the training process.

\begin{remark}
To further improve numerical accuracy, we adopt a hybrid optimization strategy. The Adam optimizer is employed in the early stage to rapidly capture the coarse structure of the solution, while the L-BFGS optimizer is used in the later stage to refine the solution and achieve improved convergence with lower residuals.
\end{remark}

\begin{algorithm}[H]
\SetAlgoLined
\DontPrintSemicolon
\caption{Deep Ritz Method with Augmented Lagrangian and Fourier Feature Mapping for the Cahn-Hilliard Equation}\label{Algorithm}
\KwIn{Computational domain $\Omega$, target mass $m_0$, maximum frequency $f_m$, outer and inner iterations $K, T$, initial penalty parameter $\mu_0$, growth rate $\rho$}
\KwOut{Trained Neural Network $u_{\theta}$}

\tcp{Step 1: Initialization}
Construct Fourier feature mapping $\Phi(x)$ with maximum frequency $f_m$\;
Initialize the neural network $u_{\theta}$ using a \textbf{GlobalResNet} architecture with \textbf{Kaiming Initialization}\;
Set Lagrange multiplier $\lambda \leftarrow 0$ and penalty parameter $\mu \leftarrow \mu_0$\;

\tcp{Step 2: Augmented Lagrangian Optimization}
\For{$k = 1, 2, \dots, K$}{
    \For{$t = 1, 2, \dots, T$}{
        Sample interior points $\mathbf{x}_i \in \Omega$ and compute $\Phi(\mathbf{x}_i)$\;
        Evaluate the Ritz energy functional: 
        $$J(u_{\theta}) = \int_{\Omega} \left( \frac{\varepsilon^2}{2} |\nabla u_{\theta}|^2 + \frac{1}{4}(u_{\theta}^2 - 1)^2 \right) d \mathbf{x}$$\;
        Compute the mass constraint: 
        $$c(\theta) = \frac{1}{N}\sum_{i=1}^N u_\theta(\mathbf{x}_i) - m_0$$\;
        Update $\theta$ by minimizing the augmented Lagrangian loss in \eqref{L_tot_FFM}: $$\mathcal{L}_{\mathrm{total}}(\theta) = J(u_{\theta}) + \lambda c(\theta) + \frac{\mu}{2} c(\theta)^2$$
        using the Adam optimizer\;
    }
    \tcp{Step 3: Multiplier and Penalty Update}
    Update the Lagrange multiplier: $\lambda \leftarrow \lambda + \mu \, c(\theta)$\;
    Update the penalty parameter: $\mu \leftarrow \min(\rho \, \mu, \mu_{max})$\;
}

\tcp{Step 4: Solution Refinement}
Refine the solution by applying the L-BFGS optimizer to further minimize $\mathcal{L}(\theta)$\;
\Return $u_{\theta}$
\end{algorithm}

\section{Numerical Experiments}\label{Numerical_results}

In this section, we evaluate the performance of the proposed method for computing steady-state solutions of the Cahn–Hilliard equation in one-, two-, and three-dimensional settings. The neural network adopts a residual architecture with a depth of six layers, each containing 100 neurons. In all the subsequent numerical results, the activation function is chosen as $\max\{\mathbf{x}^{1.5}, 0\}$, and the Adam optimizer is employed during training unless otherwise specified.

To quantitatively assess convergence, an appropriate evaluation metric is required. Since the CH equation is the gradient flow of Ginzburg-Landau free energy under $H^{-1}$ metric, a commonly used indicator is the $H^{-1}$ norm of the variational derivative (or force), denoted by $\|\Delta \delta_\phi F\|_{H^{-1}}$. In principle, this quantity vanishes as the solution approaches a steady state.
However, in the present setting, this indicator becomes unreliable. The incorporation of Fourier Feature Mappings enhances the representation of high-frequency components, while the norm of the force involves higher-order spatial derivatives, i.e., $\|\Delta \delta_\phi F\|_{H^{-1}}^2 = \inpd{-\delta_\phi F}{\Delta \delta_\phi F}_{L^2} = \| \nabla \delta_\phi F \|_{L^2}^2$. As a result, high-frequency modes are significantly amplified, and higher-order derivative terms (e.g., $\phi_{xxx}$ in one dimension) become highly sensitive. This often leads to oscillatory and unstable behavior in the force metric, which may obscure the true convergence of the solution.

To overcome this limitation, we instead adopt the difference between two successive iterations as the convergence measure:
\begin{equation}\label{Error}
\mathrm{Error} = \left( \int_\Omega \left|  u_\theta^{n+1} - u_\theta^n \right|^2 \, d\mathbf{x} \right)^{1/2}.
\end{equation}
This metric avoids the use of high-order derivatives and is therefore significantly more stable under Fourier feature representations. Consequently, it provides a more robust and reliable indicator for assessing convergence during training. Moreover, this quantity naturally approaches zero as the solution converges to a stationary state, making it consistent with the underlying variational structure.
All subsequent numerical results are evaluated using this error metric.

To further demonstrate the correctness of the numerical result by the proposed deep learning (DL) method, we also solve the one- and two- dimensional Cahn-Hilliard equation numerically by the finite difference method (FDM). 
The FDM results are taken as the reference result to demonstrate the correctness of the DL result. 
The $L^2$ error $\|\phi_{\mathrm{FDM}} - \phi_{\mathrm{DL}}\|_{L^2}$ and $L^\infty$ error $\|\phi_{\mathrm{FDM}} - \phi_{\mathrm{DL}}\|_{L^\infty}$ between the two methods are calculated for each case. 
Besides, we also calculate the free energy of the steady state by both methods and their relative error $E_{\mathrm{rel}}$ for further illustration:
$$ E_{\mathrm{rel}} = \frac{ |E_{\mathrm{FDM}} - E_{\mathrm{DL}}| }{E_{\mathrm{FDM}}}, $$
where $E_{\mathrm{FDM}}$ and $E_{\mathrm{DL}}$ are the free energy of the steady state calculated by the FDM and the DL, respectively.

In the following, we verify the proposed method by the one-, two- and three-dimensional numerical results.

\subsection{One-dimensional case}\label{subsec2}

We first consider the one-dimensional problem:
\[
F(\phi) = \int_{\Omega} \left[ \frac{\varepsilon^2}{2} |\nabla \phi|^2 + \frac{1}{4}(\phi^2-1)^2 \right] \, dx, 
\quad x \in \Omega,
\]
where $\Omega = [0,1]$ with the periodic boundary condition. The parameter is set to $\varepsilon = 0.04$, and the prescribed mass is $m = 0.6$.

\subsubsection{Trivial solution}\label{1D_numer_triv}
In the first subsection, 
we demonstrate that the trivial constant solution could be easily obtained for the one-dimensional case by the proposed framework. In this case, there is no need to apply the Fourier feature mappings technique. The periodic boundary condition is achieved by adding the boundary information as part of the loss function.
Specifically,
let $\{x_i\}_{i=1}^N \subset \Omega$ denote the interior collocation points. The total loss consists of the following components:

\begin{itemize}
    \item \textbf{Energy loss}:
    \begin{equation}\label{L_energy}
        \mathcal{L}_i = \frac{1}{N} \sum_{i=1}^N \left( 
        \frac{\varepsilon^2}{2} |\nabla u_\theta(x_i)|^2 
        + \frac{1}{4} \big(u_\theta^2(x_i) - 1\big)^2 
        \right).
    \end{equation}

    \item \textbf{Mass constraint loss}:
    \begin{equation}\label{L_mass}
        \mathcal{L}_m = \lambda \bigl( \bar{u}_\theta - m_0 \bigr) 
        + \frac{\mu}{2} \bigl( \bar{u}_\theta - m_0 \bigr)^2.
    \end{equation}

    \item \textbf{Boundary loss}:
    \[
        \mathcal{L}_b = \big| u_\theta(0) - u_\theta(1) \big|.
    \]
\end{itemize}

The total loss is defined as
\[
\mathcal{L}_{\mathrm{total}} = \mathcal{L}_i + \mathcal{L}_m + \alpha \mathcal{L}_b,
\]
where the parameters are chosen as $\lambda = 0$, $\mu = 1000$, $\alpha = 1$, with the initial augmented Lagrangian parameters $\mu_0 = 6$, $\rho = 1.2$, and $\mu_{\max} = 10$. The number of collocation points is $N = 101$. The network is trained for 500 epochs with a learning rate of $\eta = 10^{-4}$. The mean value $\bar{u}_\theta$ is approximated using Sobol quasi-Monte Carlo sampling as described in \eqref{ubar}.

\begin{figure}[H]
	\centering
	\begin{minipage}[c]{0.45\textwidth}
		\centering
		\includegraphics[width=\textwidth]{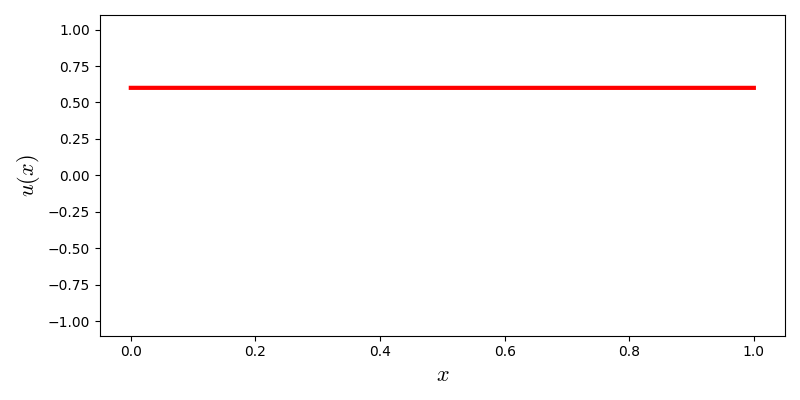}
		\subcaption{steady state}
		%\label{1d_eg1_pred}
	\end{minipage}
	\caption{Trivial steady-state solution $u \equiv 0.6$ under random initialization.}
	\label{1d_eg1}
\end{figure}

\begin{figure}[H]
	\centering
	\begin{minipage}[c]{0.30\textwidth}
		\centering
		\includegraphics[width=\textwidth]{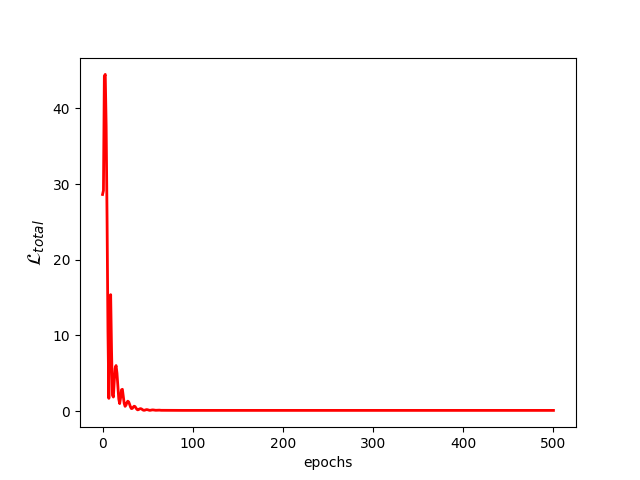}
		\subcaption{total loss}
		%\label{1d_eg1_losst}
	\end{minipage} 
    	\begin{minipage}[c]{0.30\textwidth}
		\centering
		\includegraphics[width=\textwidth]{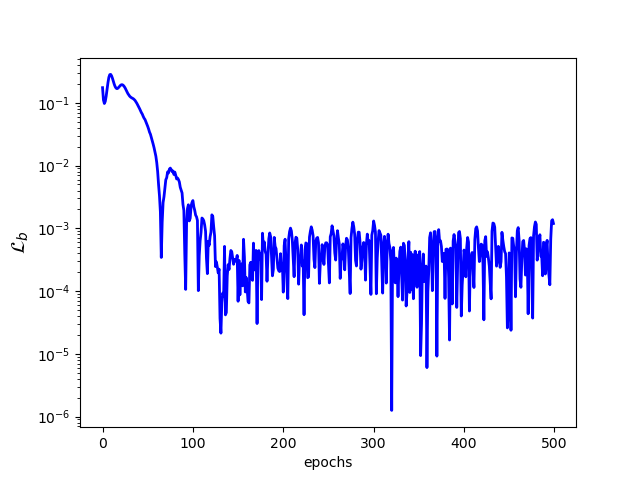}
		\subcaption{boundary loss}
		%\label{1d_eg1_lossb}
	\end{minipage} 
	\begin{minipage}[c]{0.30\textwidth}
		\centering
		\includegraphics[width=\textwidth]{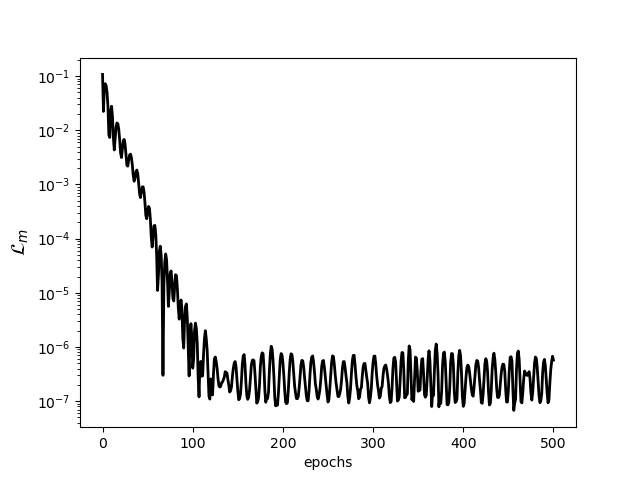}
		\subcaption{mass constraint loss}
		%\label{1d_eg1_lossm}
	\end{minipage}
	\caption{Evolution of the total loss, boundary loss, and mass constraint loss during training for the trivial solution.}
	\label{1d_eg1_losses}
\end{figure}

%Figure \ref{1d_eg1}a illustrates the random initialization of the neural network using Kaiming initialization. 
Figure \ref{1d_eg1} shows the trivial steady-state solution $u \equiv 0.6$ satisfying both the mass constraint and periodic boundary condition under random initialization using Kaiming initialization.
Figure \ref{1d_eg1_losses} presents the evolution of the total loss, boundary loss, and mass constraint loss during training. It can be observed that the total loss converges to approximately $0.1024$, which is consistent with the free energy of the steady-state solution.

\begin{comment}
\begin{figure}[H]
	\centering
	\begin{minipage}[c]{0.45\textwidth}
		\centering
		\includegraphics[width=\textwidth]{pictures/1dim/1dim_trivial_initial.png}
		\subcaption{random initializion}
		%\label{1d_eg1_initial}
	\end{minipage} 
	\begin{minipage}[c]{0.45\textwidth}
		\centering
		\includegraphics[width=\textwidth]{pictures/1dim/1dim_trivial_pred.png}
		\subcaption{steady state}
		%\label{1d_eg1_pred}
	\end{minipage}
	\caption{(a) Random initialization of the neural network using Kaiming initialization. 
(b) Recovered trivial steady-state solution $u \equiv 0.6$.}
	\label{1d_eg1}
\end{figure}    
\end{comment}

\subsubsection{Non-trivial solutions}\label{subsubsec2}

To capture the high-frequency feature of the state and get the non-trivial solution, the same technique to deal with the periodic boundary condition as the above subsection becomes insufficient. 
  To overcome this limitation, we incorporate Random Fourier Feature (RFF) mappings to enrich the representation capability of the neural network.

Since standard RFF does not inherently preserve periodicity, we adopt a simple preprocessing step by mapping the input coordinate as $x \mapsto \operatorname{mod}(x,1)$. As a result, the transformed features $\sin(2\pi \mathbf{B} x)$ and $\cos(2\pi \mathbf{B} x)$ satisfy the periodic boundary condition $u(0) = u(1)$.
In the implementation, the Fourier feature parameters are chosen as $m = 32$, $s = 3.0$, and $L_x = 1$. The outer iteration number is set to 5 and the inner iteration number to 100. The learning rate is $\eta = 10^{-4}$.  

In this setting, the total loss reduces to
\[
\mathcal{L}_{\mathrm{total}} = \mathcal{L}_i + \mathcal{L}_m,
\]
where $\mathcal{L}_i$ and $\mathcal{L}_m$ are defined the same as in \eqref{L_energy} and \eqref{L_mass}, respectively, in the previous subsection.

\begin{figure}[H]
	\centering
	\begin{minipage}[c]{0.30\textwidth}
		\centering
		\includegraphics[width=\textwidth]{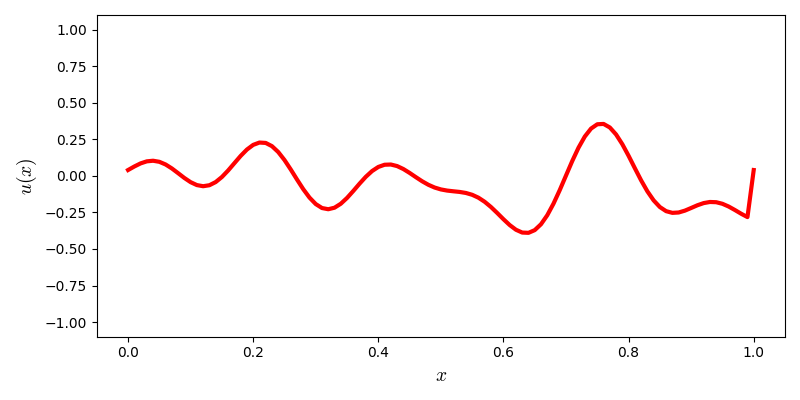}
		\subcaption{random initializion}
		%\label{1d_eg2_i}
	\end{minipage} 
	\begin{minipage}[c]{0.30\textwidth}
		\centering
		\includegraphics[width=\textwidth]{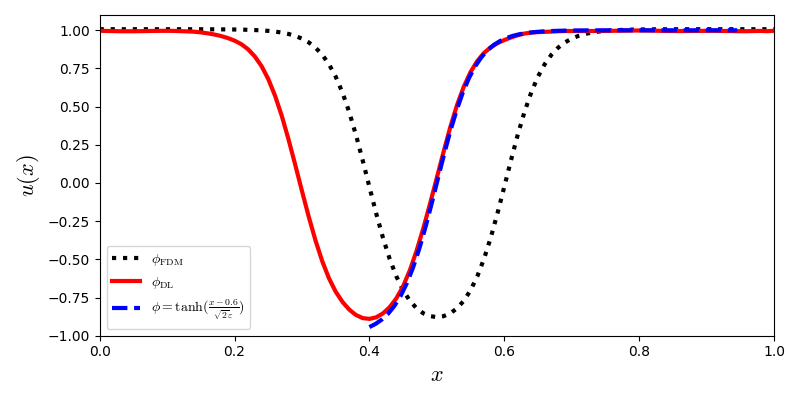}
		\subcaption{steady state}
		%\label{1d_eg2_p}
	\end{minipage}
    	\begin{minipage}[c]{0.30\textwidth}
		\centering
		\includegraphics[width=\textwidth]{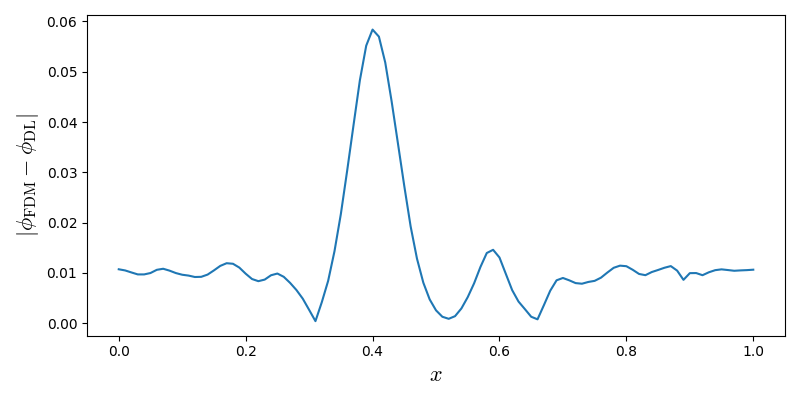}
		\subcaption{point-wise error}
		%\label{1d_eg2_p}
	\end{minipage}
	\caption{(a) Random initialization of the neural network. 
% (b) Non-trivial steady-state solution of the CH equation  enabled by Fourier feature mappings.
(b) Non-trivial steady-state solutions obtained by DL and FDM, as well as the analytical solution $\phi(x) = \tanh(\frac{x-0.6}{\sqrt{2}\epsilon})$.
(c) The point-wise absolute error $ |\phi_{\mathrm{DL}} - \phi_{\mathrm{FDM}}|$ by two methods after translational alignment.}\label{1d_eg2}
\end{figure}

Starting from the random initialization shown in Figure \ref{1d_eg2}(a), we compute the non-trivial steady-state solution, depicted as the red line in Figure \ref{1d_eg2}(b). For comparison, we also calculate the steady state using the finite difference method (FDM), shown as the black dashed line in Figure \ref{1d_eg2}(b). A phase shift is observed due to the periodic boundary condition. In the FDM computations, the number of mesh points is $N=101$ and the time step size is $\Delta t = 0.001$.
Figure \ref{1d_eg2} (c) shows the point-wise absolute error between $\phi_{\mathrm{DL}}$ and $\phi_{\mathrm{FDM}}$ after translational alignment. The absolute $L^2$ error between them is $0.0124$; see Table \ref{L2_L_infty_error}. To validate the correctness of the interface thickness for the given $\varepsilon$, we plot the analytical tanh-type structure of the Cahn-Hilliard equation 
for the one-dimensional case for comparison: 
$$\phi(x) = \tanh\left(\frac{x-0.6}{\sqrt{2}\varepsilon}\right), $$
shown as the blue dashed line in Figure \ref{1d_eg2}(b). Good agreement is observed with the numerical solution at the interface, confirming the correct interface thickness.
The corresponding free energy of the non-trivial steady state is approximately $0.074945$, with a relative error of $10^{-5}$ compared to the energy obtained by FDM; see Table \ref{Relative_erro_egy}.

\begin{figure}[H]
	\centering
	\begin{minipage}[c]{0.30\textwidth}
		\centering
		\includegraphics[width=\textwidth]{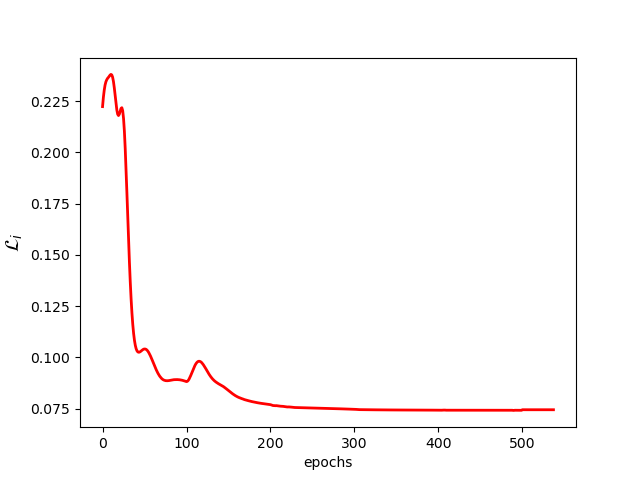}
		\subcaption{energy loss}
		%\label{1d_eg2_lossi}
	\end{minipage} 
	\begin{minipage}[c]{0.30\textwidth}
		\centering
		\includegraphics[width=\textwidth]{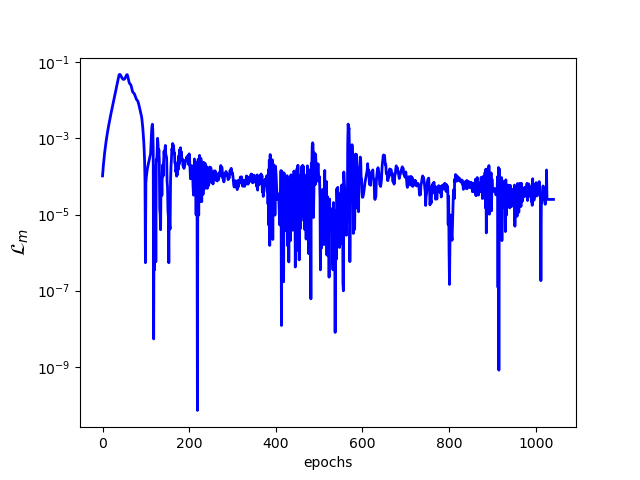}
		\subcaption{mass constraint loss}
		%\label{1d_eg2_lossm}
	\end{minipage}
	\begin{minipage}[c]{0.30\textwidth}
		\centering
		\includegraphics[width=\textwidth]{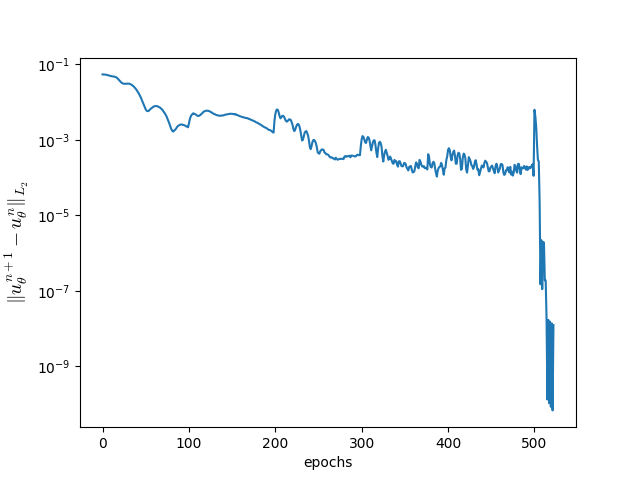}
		\subcaption{error decay}
		%\label{force}
	\end{minipage}
	\caption{Evolution of (a) energy loss, (b) mass constraint loss, and (c) error, demonstrating stable convergence to a non-trivial steady state.}
	\label{1d_eg2_losses}
\end{figure}

Figures \ref{1d_eg2_losses}(a) and (b) show the evolution of the energy loss and mass constraint loss, respectively. The total number of epochs corresponds to the product of the inner and outer iteration numbers. 
The total computational time is approximately $13.59$ seconds, demonstrating the efficiency of the proposed method. To further assess the accuracy, we plot the error defined in \eqref{Error} against the training epochs, as shown in Figure \ref{1d_eg2_losses}(c). It should be noted that, to refine the solution, the L-BFGS optimization method is employed in the later stage of the algorithm to achieve a more significant reduction in the final force values. The error exhibits a clear decay trend, indicating convergence to a steady state.

These results demonstrate that the incorporation of Fourier feature mappings effectively overcomes the spectral bias of neural networks and enables the accurate computation of non-trivial steady states.

\subsection{Two-dimensional case}\label{subsec3}

We next consider the two-dimensional case. The energy functional is given by
\begin{equation}\label{GL_egy_2D}
    F(\phi) = \int_{\Omega} \left[ \frac{\varepsilon^2}{2} |\nabla \phi|^2 + \frac{1}{4}(\phi^2-1)^2 \right] \, d \mathbf{x}, 
\quad \mathbf{x} \in \Omega,    
\end{equation}
where $\Omega = [0,1]\times [0,1]$ and $\varepsilon = 0.01$. The prescribed mass is $m_0 = 0.02$. Periodic boundary conditions are imposed in both spatial directions:
\[
\phi(0,y) = \phi(1,y), \quad y \in [0,1], \qquad 
\phi(x,0) = \phi(x,1), \quad x \in [0,1].
\]
As in the one-dimensional case, both trivial and non-trivial steady-state solutions may exist.

\subsubsection{Trivial solution} \label{2D_numer_triv}

We first consider the trivial solution. The periodic boundary conditions are enforced directly, following the same strategy as in the one-dimensional case. The total loss function consists of three components:

\begin{itemize}
    \item \textbf{Energy loss}:
    \[
        \mathcal{L}_i = \frac{1}{N^2}\sum_{i=1}^N \sum_{j=1}^N  
        \left( 
        \frac{\varepsilon^2}{2} |\nabla u_\theta(x_i,y_j)|^2 
        + \frac{1}{4} \big(u_\theta^2(x_i,y_j) - 1\big)^2 
        \right).
    \]

    \item \textbf{Mass constraint loss}:
    \[
       \mathcal{L}_m = \lambda \bigl( \bar{u}_\theta - m_0 \bigr) 
       + \frac{\mu}{2} \bigl( \bar{u}_\theta - m_0 \bigr)^2.
    \]

    \item \textbf{Boundary loss}:
    \[
        \mathcal{L}_b = \frac{1}{N} \sum_{i=1}^N \big(u_\theta(x_i,0) - u_\theta(x_i,1)\big)^2 
        + \frac{1}{N} \sum_{j=1}^N \big(u_\theta(0,y_j) - u_\theta(1,y_j)\big)^2.
    \]
\end{itemize}

The total loss is defined as
\[
    \mathcal{L}_{\mathrm{total}} = \mathcal{L}_i + \mathcal{L}_m + \alpha \mathcal{L}_b,
\]
where the parameters are chosen as $\lambda = 0$, $\mu = 20$, $\mu_0 = 1$, $\rho = 1.2$, $\mu_{\max} = 2$, and $\alpha = 10$. 

We employ $N = 202$ uniformly distributed collocation points in each spatial direction. The numbers of outer and inner iterations are set to 5 and 100, respectively, with a learning rate of $\eta = 10^{-3}$.

\begin{figure}[H]
	\centering
	\begin{minipage}[c]{0.30\textwidth}
		\centering
		\includegraphics[width=\textwidth]{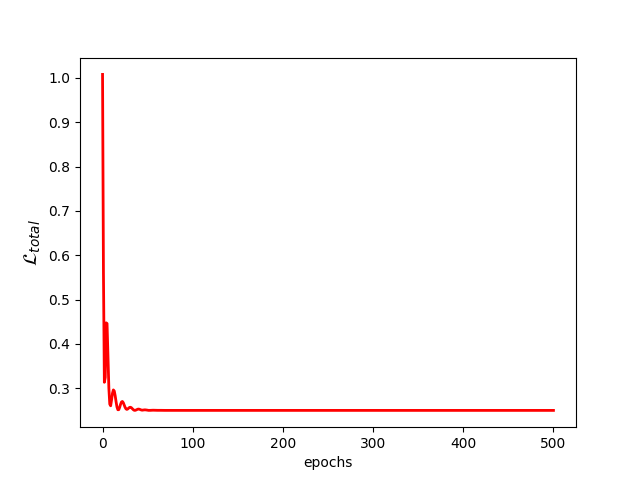}
		\subcaption{total loss}
		%\label{fig_E2_1}
	\end{minipage} 
    	\begin{minipage}[c]{0.30\textwidth}
		\centering
		\includegraphics[width=\textwidth]{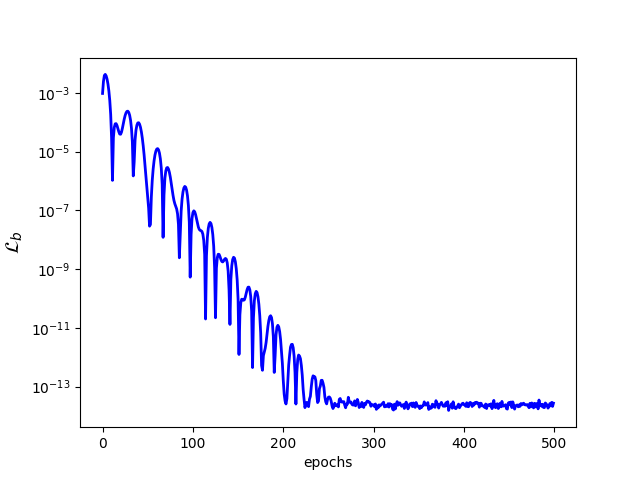}
		\subcaption{boundary loss}
		%\label{fig_E2_2}
	\end{minipage} 
	\begin{minipage}[c]{0.30\textwidth}
		\centering
		\includegraphics[width=\textwidth]{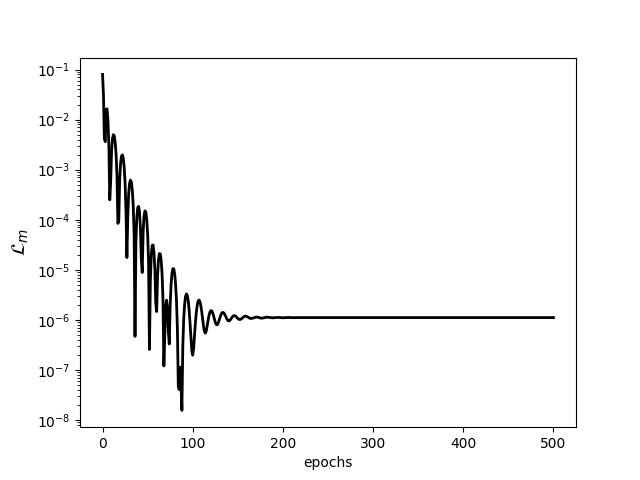}
		\subcaption{mass constraint loss}
		%\label{fig_E2_3}
	\end{minipage}
	\caption{Evolution of the total loss, boundary loss, and mass constraint loss during training.}
	\label{2d_eg1_losses}
\end{figure}

Starting from a random initialization, the network converges to the constant steady-state solution $u \equiv 0.02$, which satisfies both the mass constraint and periodic boundary conditions. This is consistent with the FDM result, where the free energy is $0.2498$.
Figure~\ref{2d_eg1_losses} illustrates the evolution of the total loss, boundary loss, and mass constraint loss during training, all of which exhibit stable decay behavior. 
However, as in the one-dimensional case, this formulation tends to favor trivial constant solutions and fails to capture more complex non-trivial patterns.

% For the Lagrangian Method,  in example 2 and 3, the outer iterations is 12 while the others are 10. And inner iterations is 100 for all. $\mu_0=1, \rho=1.2, \mu_{max}=2$. For the Fourier feature mapping, the maximum frequency is 3, $L_x=L_y=1$. The learning rate is 0.001.
% % In the two-dimensional case, we found four steady states: constant solution(ex.1), elliptic solution $\phi=1$ (ex.2, ex.5), elliptic solution $\phi=-1$ (ex.3, ex.4), and striped solution(ex.6). 
% Example 3, 4 and 5 cost about 60 seconds to get the results.

%\begin{figure}[H]
%    \centering
%    \includegraphics[width=0.3\linewidth]{pictures/2dim/2dim_trivial_initial.png}
%    \caption{random initializion}
%    \label{2d_eg1}
%\end{figure}

%Fig.\ref{2d_eg1} shows using kaiming to randomly initialize the neural network, and the final steady-state solution is 0.02. 

\subsubsection{Non-trivial solutions}

To compute non-trivial steady-state solutions in the two-dimensional case, we employ Fourier feature mappings to enhance the expressive capability of the neural network.

In the one-dimensional case, non-trivial solutions are obtained using standard random Fourier features (RFF), which provide flexible bandwidth control and improved representation of high-frequency structures. However, in higher-dimensional settings, to better preserve periodicity and separability along each coordinate direction, we instead adopt classical Fourier feature mappings based on integer-frequency Fourier series.

Specifically, the maximum frequency is chosen as $f_m = 3$, with domain lengths $L_x = L_y = 1$. Under this construction, periodic boundary conditions are naturally satisfied, and the total loss function reduces to
\[
\mathcal{L}_{\mathrm{total}} = \mathcal{L}_i + \mathcal{L}_m,
\]
where $\mathcal{L}_i$ and $\mathcal{L}_m$ are defined as in Subsection~\ref{2D_numer_triv}.

In the deep ritz method, the numbers of outer and inner iterations are set to 10 and 100, respectively, with a learning rate of $\eta = 10^{-3}$. Below, we present several representative numerical examples of non-trivial steady states. For all cases, the method converges within approximately $60$ seconds, demonstrating its computational efficiency. Additionally, we compare our results with those obtained from the finite difference method (FDM). In the FDM computations, we use $Nx = 202, Ny=202$ mesh points and a time step size of $\Delta t = 10^{-2}$ for all the cases.

\paragraph{\textit{Case 1}} \label{2D_ex2}

\begin{figure}[H]
	\centering
	\begin{minipage}[c]{0.23\textwidth}
		\centering
		\includegraphics[width=\textwidth]{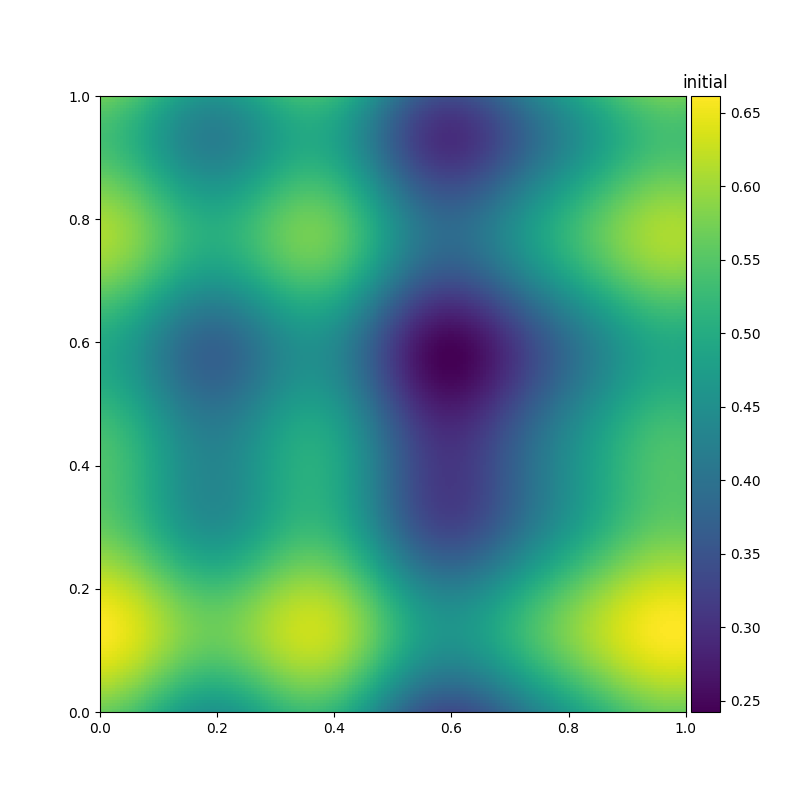}
		\subcaption{random initialization}
		\label{2d_eg1_initial}
	\end{minipage} 
	\begin{minipage}[c]{0.23\textwidth}
		\centering
		\includegraphics[width=\textwidth]{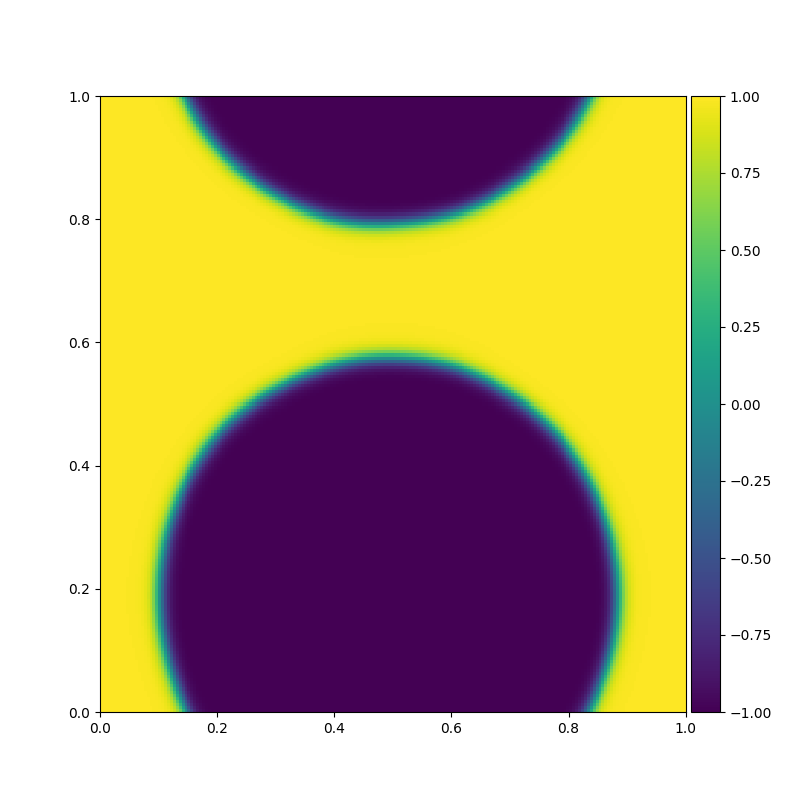}
		\subcaption{steady state(DL)}
		\label{2g_eg1_pred}
	\end{minipage}
	\begin{minipage}[c]{0.23\textwidth}
		\centering
		\includegraphics[width=\textwidth]{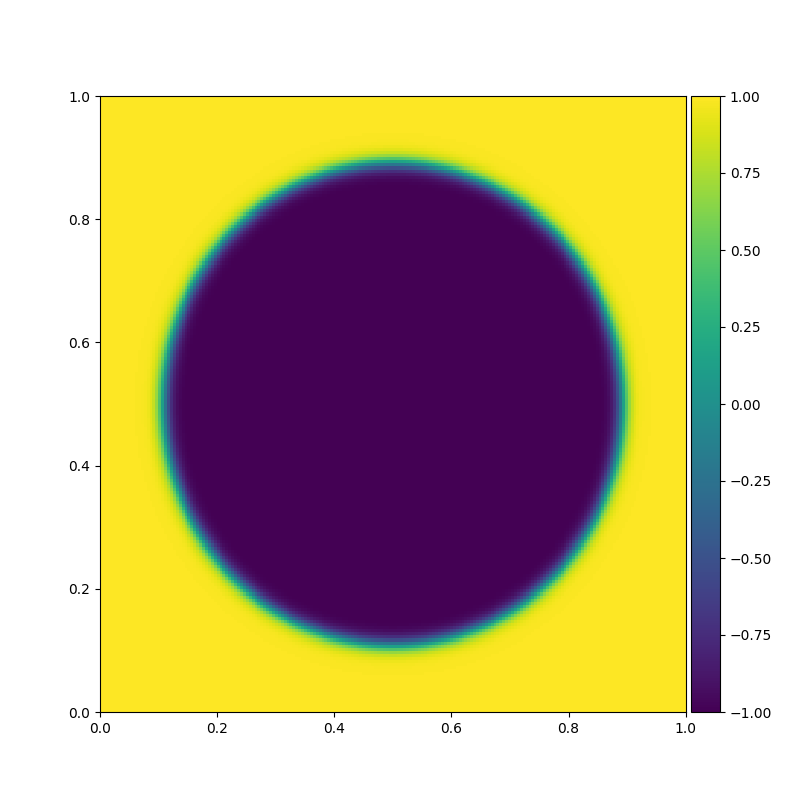}
		\subcaption{steady state(FDM)}
		\label{2g_eg1_error}
	\end{minipage}    
	\begin{minipage}[c]{0.23\textwidth}
		\centering
		\includegraphics[width=\textwidth]{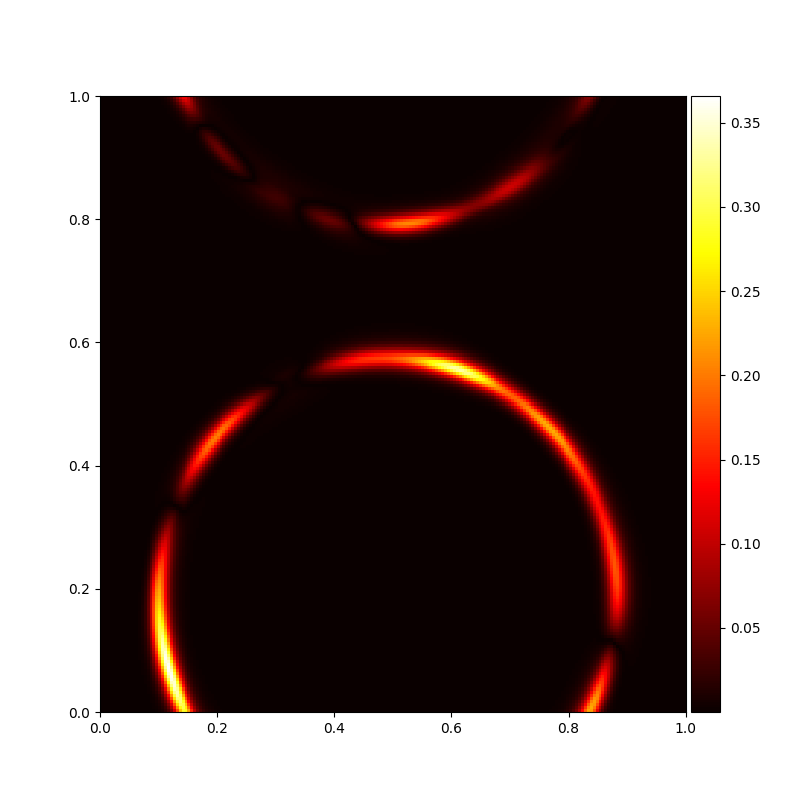}
		\subcaption{point-wise error}
		\label{2g_eg1_error}
	\end{minipage}    
	\caption{Case 1 (2D): (a) Random initialization of the neural network. (b) Non-trivial droplet-type steady state of the Cahn–Hilliard equation enabled by Fourier feature mappings. (c) Reference steady state solution from FDM. (d) Translational-alignment-corrected point-wise absolute error between (b) and (c).}\label{2d_eg4}
\end{figure}

Figure~\ref{2d_eg4}(b) shows a representative non-trivial steady-state solution obtained from the random initialization displayed in Figure~\ref{2d_eg4}(a). The solution exhibits a typical phase-separated configuration, consisting of a region where $\phi\approx -1$ embedded within a surrounding region where $\phi\approx 1$, with the two phases connected by a thin and smooth interfacial layer. 
This phase profile is consistent with the numerical solution calculated by FDM, shown as in Figure \ref{2d_eg4}(c) after a translation. Figure \ref{2d_eg4}(d) presents the point-wise absolute error between the two phase profiles from DL and FDM after translational alignment. The $L^2$ norm error between them is $0.0083$; see Table \ref{L2_L_infty_error}.
The free energy of the steady state is approximately $0.0232$, with a relative error of $4.33 \times 10^{-3}$ with respect to the steady state obtained via FDM (see Table \ref{Relative_erro_egy}).

\begin{figure}[H]
	\centering
	\begin{minipage}[c]{0.30\textwidth}
		\centering
		\includegraphics[width=\textwidth]{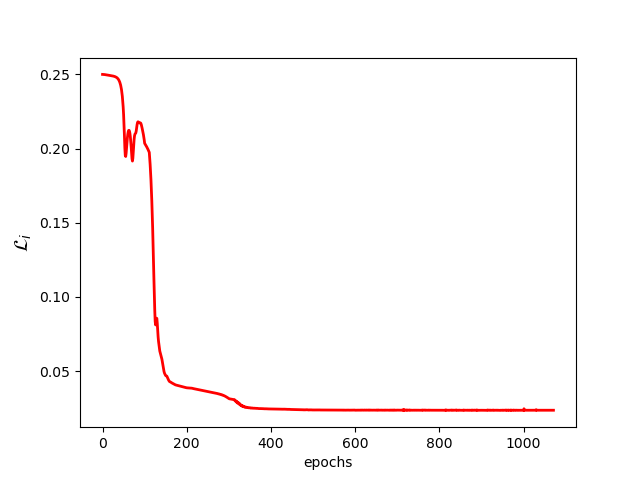}
		\subcaption{energy loss}
		\label{2d_eg1_lossi}
	\end{minipage} 
	\begin{minipage}[c]{0.30\textwidth}
		\centering
		\includegraphics[width=\textwidth]{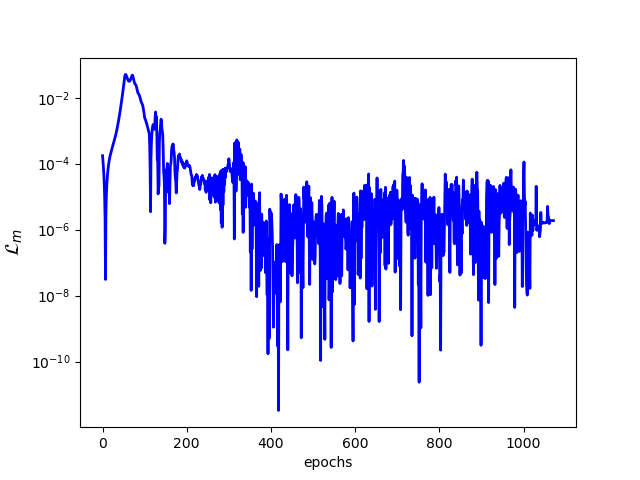}
		\subcaption{mass constraint loss}
		\label{2d_eg1_lossm}
	\end{minipage}
	\begin{minipage}[c]{0.30\textwidth}
		\centering
		\includegraphics[width=\textwidth]{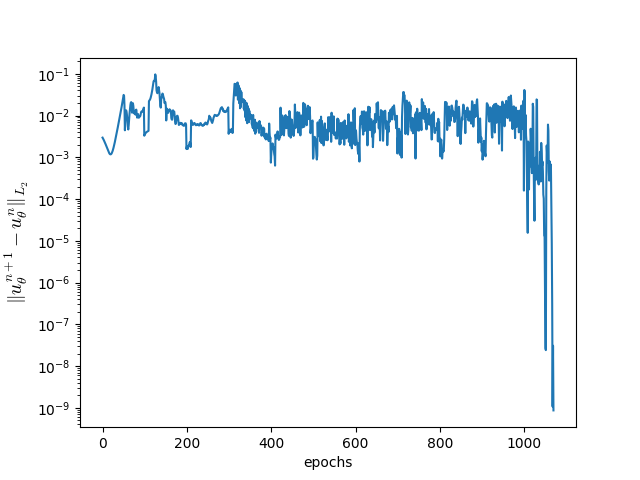}
		\subcaption{error decay}
		\label{2d_eg1_force}
	\end{minipage}    
	\caption{Case 1 (2D): Evolution of (a) energy loss, (b) mass constraint loss and (c) error during training for the droplet-type non-trivial solution.}
	\label{2d_eg4_losses}
\end{figure}

Figures~\ref{2d_eg4_losses}(a) and (b) illustrate the evolution of the energy loss and mass constraint loss, respectively. The energy loss decreases and converges to a stable value, and the mass constraint loss decreases with an average magnitude on the order of $10^{-7}$. Figure~\ref{2d_eg4_losses}(c) presents the decay of the error defined in \eqref{Error}, further confirming the convergence and precision of the calculated solution.
These results demonstrate that the proposed method is capable of efficiently capturing complex non-trivial steady states in higher-dimensional settings.

%%%%%%%%%%%%%%%%%%%%%%%%%%%%%%%%%%%%%%%%%%%%%%%%%%%%%%

\paragraph{\textit{Case 2}} \label{2D_ex3}

\begin{figure}[H]
	\centering
	\begin{minipage}[c]{0.23\textwidth}
		\centering
		\includegraphics[width=\textwidth]{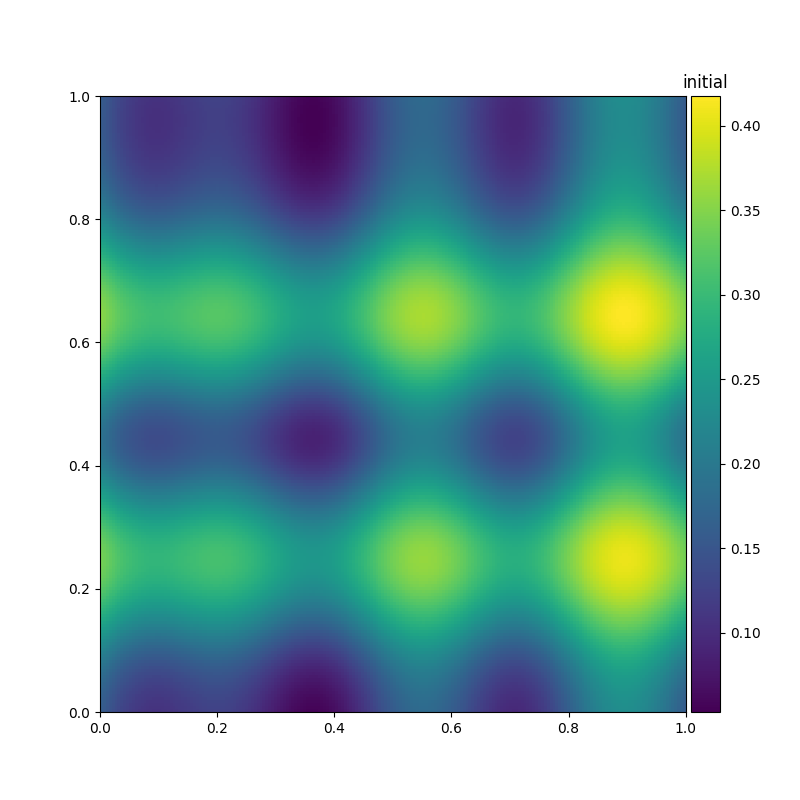}
		\subcaption{random initializion}
		\label{2d_eg2_initial}
	\end{minipage} 
	\begin{minipage}[c]{0.23\textwidth}
		\centering
		\includegraphics[width=\textwidth]{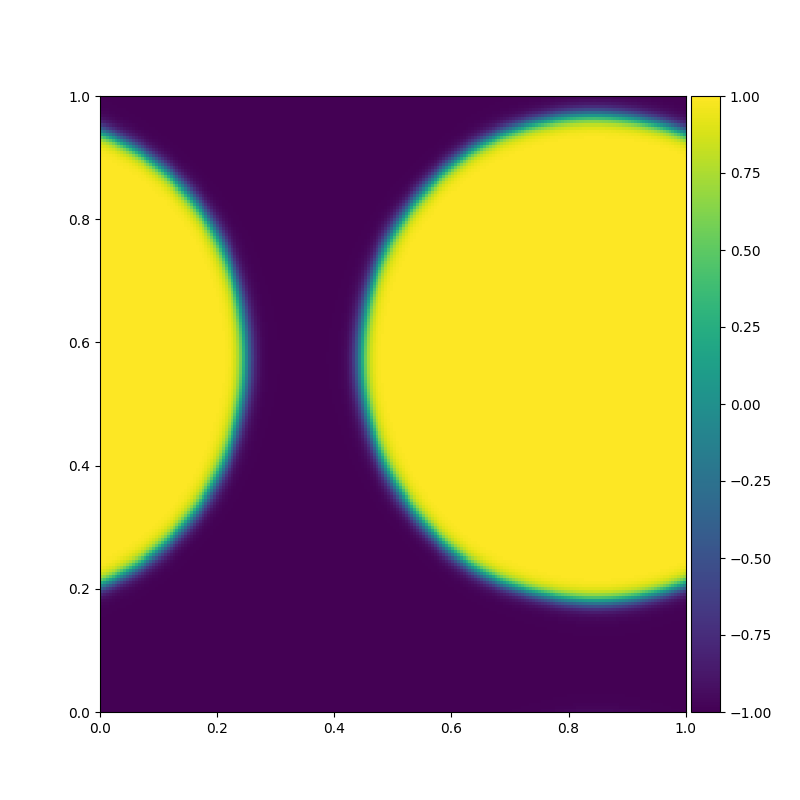}
		\subcaption{steady state(DL)}
		\label{fig_droplet_2}
	\end{minipage}
	\begin{minipage}[c]{0.23\textwidth}
		\centering
		\includegraphics[width=\textwidth]{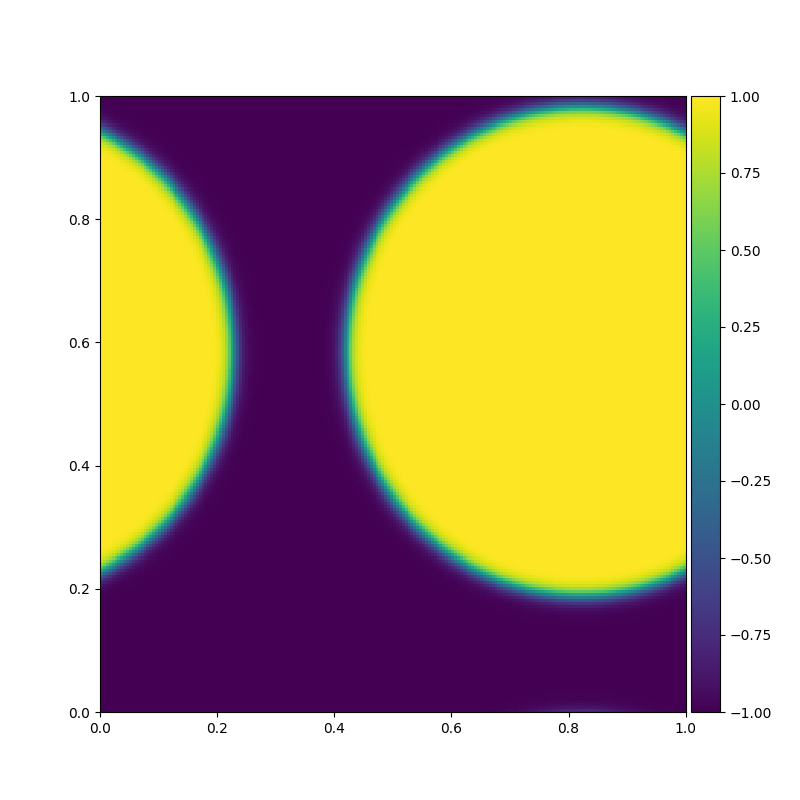}
		\subcaption{steady state(FDM)}
		%\label{fig_fdm_2}
	\end{minipage}
	\begin{minipage}[c]{0.23\textwidth}
		\centering
		\includegraphics[width=\textwidth]{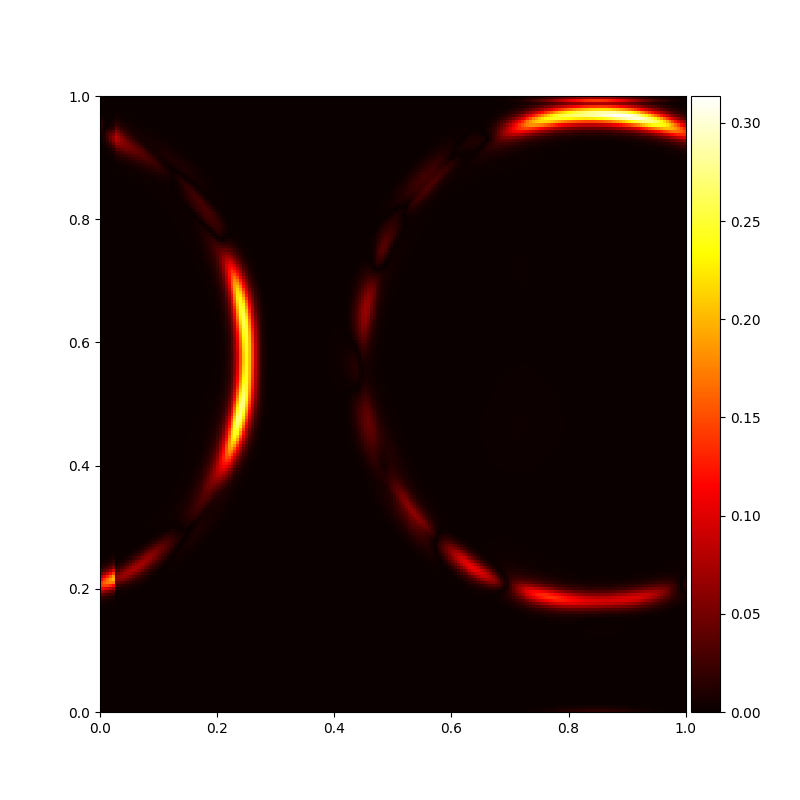}
		\subcaption{point-wise error}
		\label{fig_error_2}
	\end{minipage}
	\caption{Case 2 (2D): (a) Random initialization of the neural network. (b) Inverted droplet-type steady-state solution enabled by Fourier feature mappings, where the two phases are exchanged compared to Case~1. (c) Reference steady state solution from FDM. (d) Point-wise absolute error between (b) and (c). }
	\label{2d_eg5}
\end{figure}

Figure~\ref{2d_eg5}(b) presents another non-trivial steady-state solution obtained from the random initialization shown in Figure~\ref{2d_eg5}(a). Similarly to Case~1, the solution exhibits a phase-separated structure characterized by a droplet-like morphology. However, in contrast to Case~1, the roles of the two phases are reversed: the inner region and the surrounding phase exchange signs, resulting in an inverted droplet configuration. This is consistent with the FDM result, shown in Figure \ref{2d_eg5} (c). Figure \ref{2d_eg5} (d) shows the point-wise absolute error between the phase profiles in (b) and (c). The $L^2$ norm error between them is $0.0073$; see Table \ref{L2_L_infty_error}.
The free energy of this steady state is approximately $0.0236$, with a relative error of $4.26\times 10^{-3}$ with respect to the steady state obtained by FDM (see Table \ref{Relative_erro_egy}).

\begin{figure}[H]
	\centering
	\begin{minipage}[c]{0.30\textwidth}
		\centering
		\includegraphics[width=\textwidth]{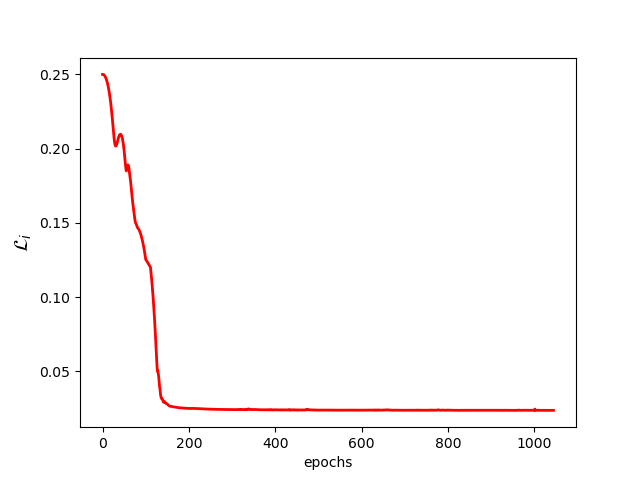}
		\subcaption{energy loss}
		\label{2d_eg2_lossi}
	\end{minipage} 
	\begin{minipage}[c]{0.30\textwidth}
		\centering
		\includegraphics[width=\textwidth]{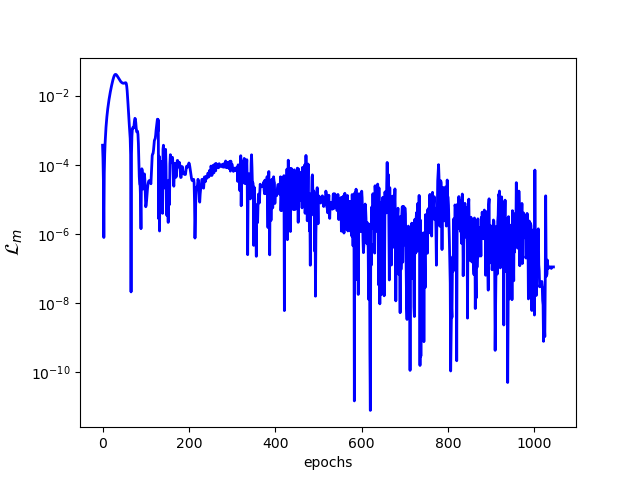}
		\subcaption{mass constraint loss}
		\label{2d_eg2_lossm}
	\end{minipage}
	\begin{minipage}[c]{0.30\textwidth}
		\centering
		\includegraphics[width=\textwidth]{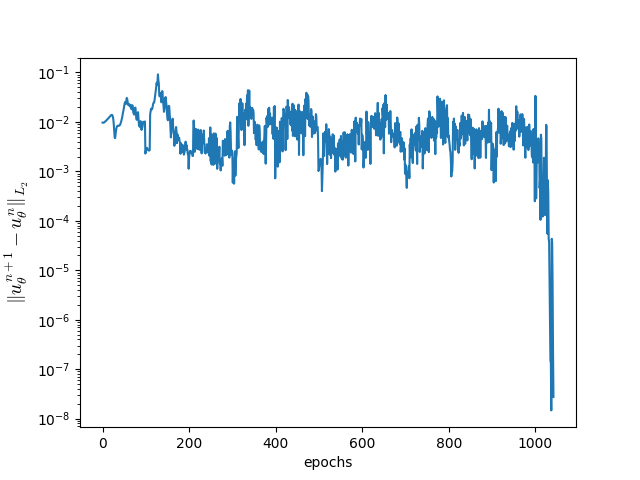}
		\subcaption{error decay}
		\label{2d_eg2_force}
	\end{minipage}
	\caption{Case 2 (2D): Evolution of (a) energy loss, (b) mass constraint loss, and (c) error during training for the inverted droplet-type solution.}
	\label{2d_eg5_losses}
\end{figure}

Figures~\ref{2d_eg5_losses}(a) and (b) show the evolution of the energy loss and mass constraint loss, respectively, both exhibiting stable convergence. The error decay, shown in Figure~\ref{2d_eg5_losses}(c), further confirms the accuracy and convergence of the solution.

%%%%%%%%%%%%%%%%%%%%%%%%%%%%%%%%%%%%%%%%%%%%%%%%%%%%%%

\paragraph{\textit{Case 3}} \label{2D_ex4}

\begin{figure}[H]
	\centering
	\begin{minipage}[c]{0.23\textwidth}
		\centering
		\includegraphics[width=\textwidth]{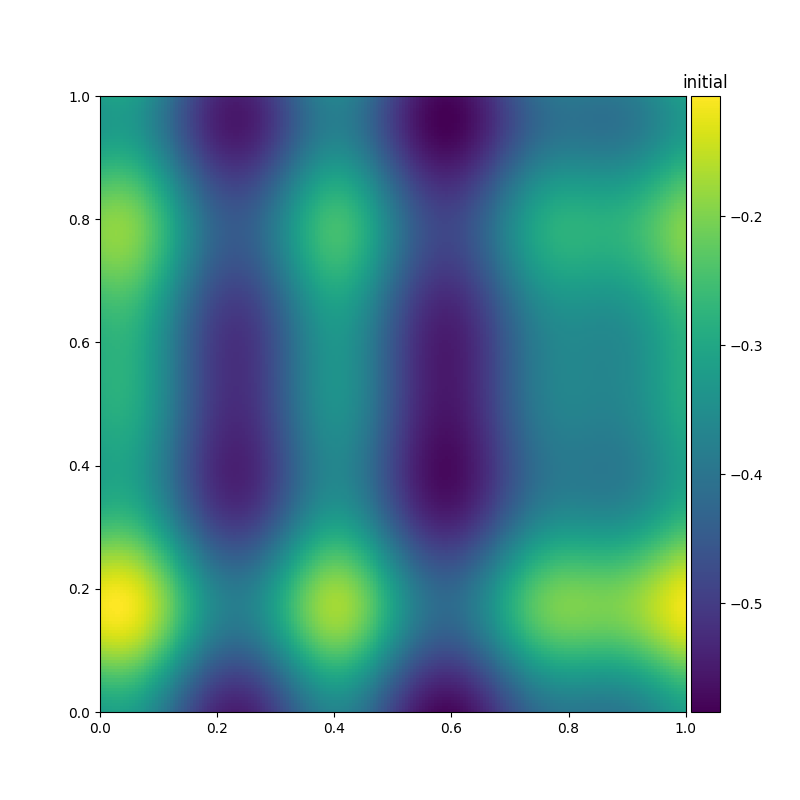}
		\subcaption{random initializion}
		\label{2d_eg3_initial}
	\end{minipage} 
	\begin{minipage}[c]{0.23\textwidth}
		\centering
		\includegraphics[width=\textwidth]{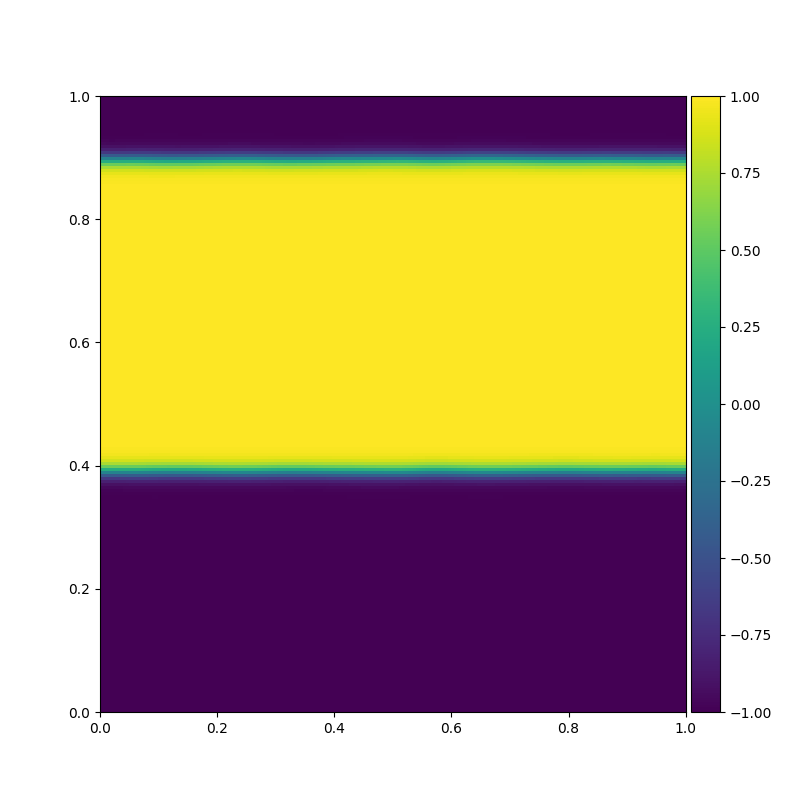}
		\subcaption{steady state(DL)}
		\label{fig_2d_strip}
	\end{minipage}
	\begin{minipage}[c]{0.23\textwidth}
		\centering
		\includegraphics[width=\textwidth]{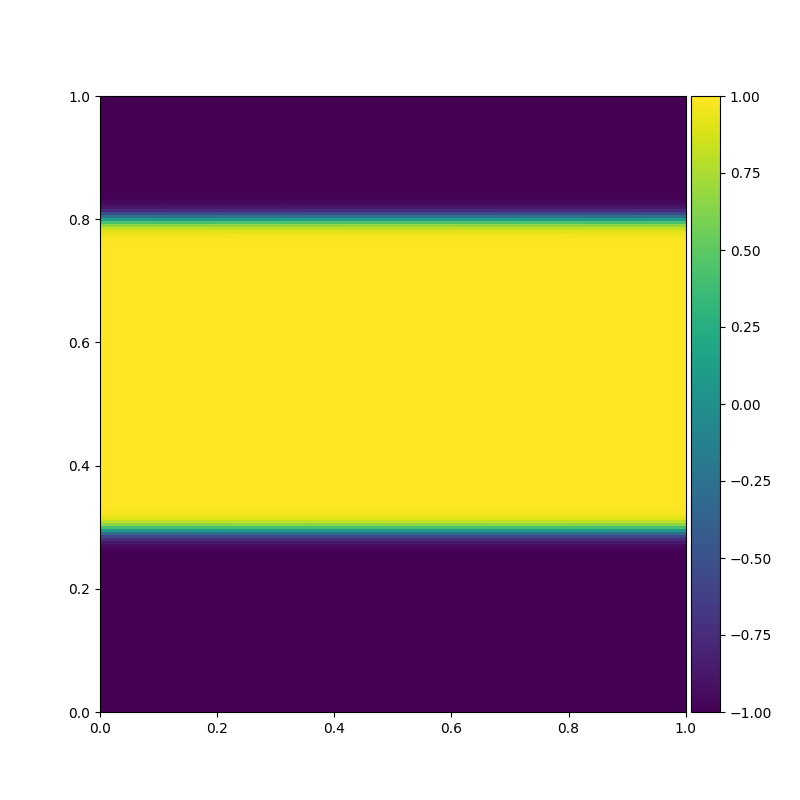}
		\subcaption{steady state(FDM)}
		\label{fig_2d_strip_fdm}
	\end{minipage}
	\begin{minipage}[c]{0.23\textwidth}
		\centering
		\includegraphics[width=\textwidth]{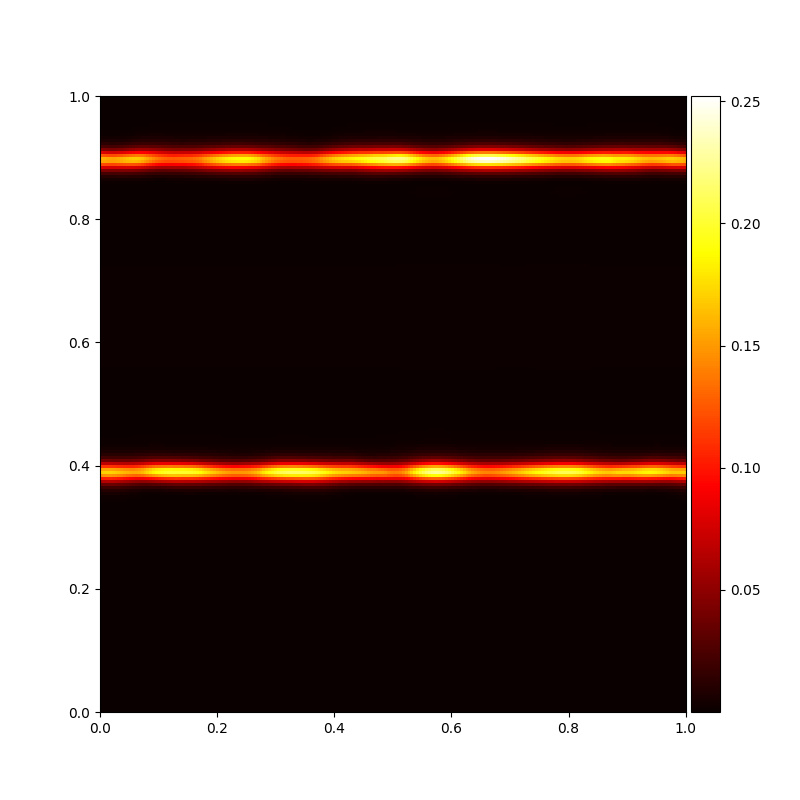}
		\subcaption{absolute error}
		\label{fig_2d_strip_error}
	\end{minipage}
	\caption{Case 3 (2D): (a) Random initialization of the neural network. (b) Lamellar (striped) steady-state solution enabled by Fourier feature mappings, exhibiting alternating layered phase separation. (c) Reference steady state from FDM. (d) Translational-alignment-corrected point-wise absolute error between (b) and (c). }
	\label{2d_eg6}
\end{figure}

From a different random initialization (Figure~\ref{2d_eg6}(a)), the method converges to a qualitatively distinct non-trivial steady state, shown in Figure \ref{2d_eg6}(b). Unlike the droplet-type structures observed in Cases~1 and~2, this solution exhibits a lamellar (striped) phase separation pattern, characterized by alternating layers of the two phases across the domain.
Such lamellar structures are commonly observed in phase separation processes and correspond to a different class of local energy minimizers. This profile is consistent with the phase profile calculated by FDM, shown in Figure \ref{2d_eg6}(c). Figure \ref{2d_eg6}(d) shows the point-wise absolute error between the two profiles in (b) and (c) after translational alignment. The $L^2$ norm error between them is $0.0105$; see Table \ref{L2_L_infty_error}. 
The free energy is approximately $0.0187$; with a relative error of $1.08 \times 10^{-2}$ with respect to the steady state obtained by FDM; see Table \ref{Relative_erro_egy}.

\begin{figure}[H]
	\centering
	\begin{minipage}[c]{0.30\textwidth}
		\centering
		\includegraphics[width=\textwidth]{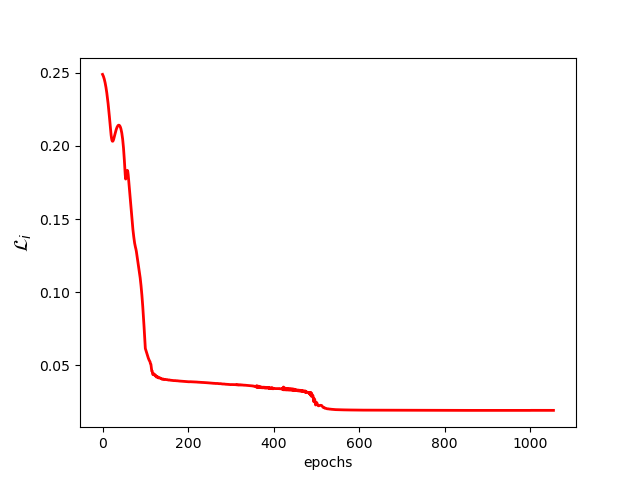}
		\subcaption{energy loss}
		\label{2d_eg3_lossi}
	\end{minipage} 
	\begin{minipage}[c]{0.30\textwidth}
		\centering
		\includegraphics[width=\textwidth]{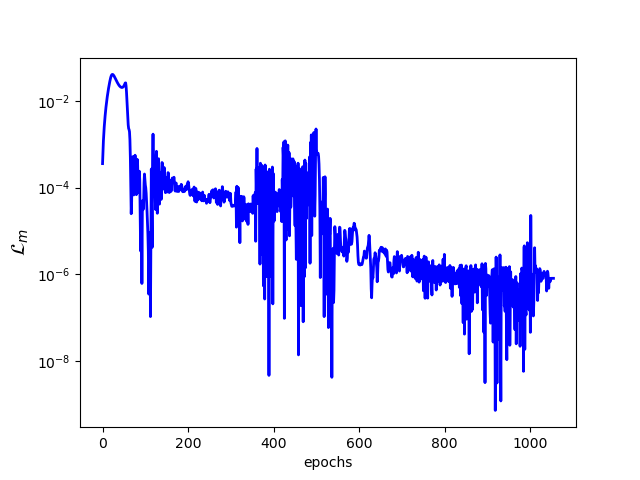}
		\subcaption{mass constraint loss}
		\label{2d_eg3_lossm}
	\end{minipage}
    	\begin{minipage}[c]{0.30\textwidth}
		\centering
		\includegraphics[width=\textwidth]{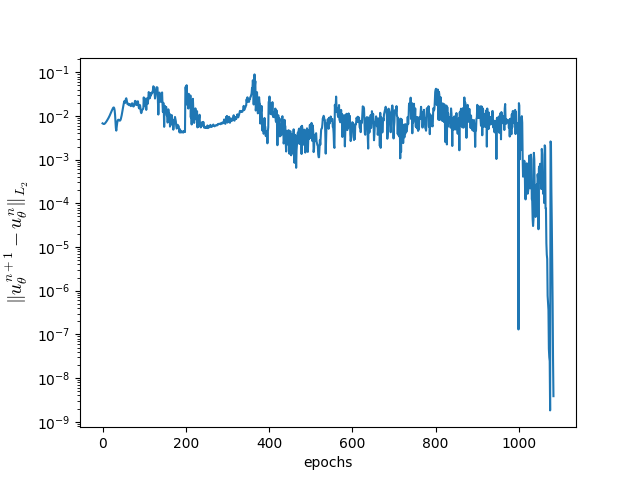}
		\subcaption{error decay}
		\label{2d_eg3_force}
	\end{minipage}
	\caption{Case 3 (2D): Evolution of (a) energy loss, (b) mass constraint loss, and (c) error during training for the lamellar steady-state solution.}
	\label{2d_eg6_losses}
\end{figure}

Figures~\ref{2d_eg6_losses}(a) and (b) illustrate the decay of the energy loss and mass constraint loss, respectively, both showing stable convergence behavior. The error curve in Figure~\ref{2d_eg6_losses}(c) further verifies the accuracy of the computed solution.

\begin{table}[H]
\centering
\caption{Comparison between the phase profiles by DL and FDM.}\label{L2_L_infty_error}
\begin{tabular}{l|l|l|l|l}
$\| \phi_{\mathrm{DL}} - \phi_{\mathrm{FDM}}\|$   & 1D(nontrivial & 2D(case1) & 2D(case2) & 2D(case3) \\ \hline
$L_2$ error                        & 0.0124        & 0.0083& 0.0073    & 0.0105    \\ \hline
$L_{\infty}$ error & 0.0584        & 0.3659    & 0.3137    & 0.2521   
\end{tabular}
\end{table}

\begin{table}[H]
\centering
\caption{Energies computed by Deep Ritz method and FDM and their relative error.}\label{Relative_erro_egy}
\begin{tabular}{c|l|l|l|l}
Energy   & 1D(nontrivial)  & 2D(case1)  & 2D(case2)  & 2D(case3) \\ \hline
Deep Ritz      & 0.074945  & 0.0232   & 0.0236    & 0.0187   \\ \hline
FDM   & 0.074947      & 0.0231   & 0.0235    & 0.0185       \\ \hline
Relative error &\num{2.67e-5}& \num{4.33e-3} & \num{4.26e-3} & \num{1.08e-2}
\end{tabular}
\end{table}

\subsection{Three-dimensional case}\label{subsec4}

We next consider the three-dimensional Cahn--Hilliard problem. The energy functional remains the same as in \eqref{GL_egy_2D}, with the computational domain $\Omega = [0,1]^3$, parameter $\varepsilon = 0.01$, and prescribed mass $m_0 = 0.02$. The periodic boundary conditions are given by
\begin{align*}
\phi(0,y,z) &= \phi(1,y,z), \quad (y,z)\in [0,1]^2, \\
\phi(x,0,z) &= \phi(x,1,z), \quad (x,z)\in [0,1]^2, \\
\phi(x,y,0) &= \phi(x,y,1), \quad (x,y)\in [0,1]^2.
\end{align*}

As in the one- and two-dimensional cases, a trivial constant solution $u \equiv 0.02$ exists. In this subsection, we focus on computing non-trivial steady-state solutions using the separable Fourier feature mapping technique.
With the incorporation of Fourier features, the periodic boundary conditions are naturally satisfied. Therefore, the total loss functional consists only of the energy term and the mass constraint term:
\[
\mathcal{L}_{\mathrm{total}} = \mathcal{L}_i + \mathcal{L}_m,
\]
where
\begin{itemize}
    \item \textbf{Energy loss}:
    \[
    \mathcal{L}_i = \frac{1}{N}\sum_{i=1}^N 
    \left( 
    \frac{\varepsilon^2}{2} |\nabla u_\theta(\mathbf{x}_i)|^2 
    + \frac{1}{4} \big(u_\theta^2(\mathbf{x}_i) - 1\big)^2 
    \right),
    \]

    \item \textbf{Mass constraint loss}:
    \[
    \mathcal{L}_m = \lambda \big( \bar{u}_\theta - m_0 \big) 
    + \frac{\mu}{2} \big( \bar{u}_\theta - m_0 \big)^2.
    \]
\end{itemize}
Due to the rapid growth in the number of grid points in three dimensions, uniform discretization becomes computationally prohibitive. To address this issue, we employ Sobol low-discrepancy sequences to generate collocation points, which provide high-quality space-filling sampling with relatively small sample sizes. In our experiments, we use $N = 100{,}000$ sampling points.

For the separable Fourier feature mapping, the maximum frequency is set to $f_m = 3$, with domain lengths $L_x = L_y = L_z = 1$. The augmented Lagrangian parameters are chosen as $\mu_0 = 1$, $\rho = 1.2$, and $\mu_{\max} = 2$. The optimization is performed with 10 outer iterations and 100 inner iterations, and the learning rate is set to $\eta = 5 \times 10^{-4}$.

Starting from different random initializations, the proposed method is able to converge to multiple distinct non-trivial steady states in three dimensions (see the following figures), demonstrating its capability to capture complex phase separation structures in high-dimensional settings.

%%%%%%%%%%%%%%%%%%%%%%%%%%%%%%%%%%%%%%%%%%%%%%%%%%%%%%
\paragraph{\textit{Case 1}}

\begin{figure}[H]
	\centering
	\begin{minipage}[c]{0.30\textwidth}
		\centering
		\includegraphics[width=\textwidth]{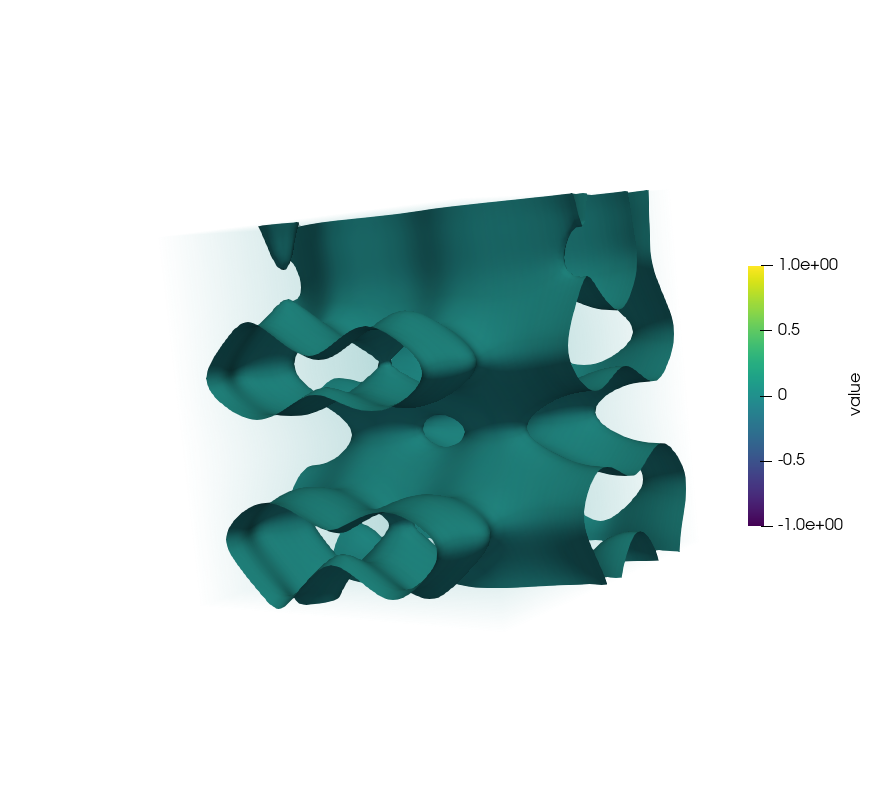}
		\subcaption{random initializion}
		\label{3d_eg1_initial}
	\end{minipage} 
	\begin{minipage}[c]{0.30\textwidth}
		\centering
		\includegraphics[width=\textwidth]{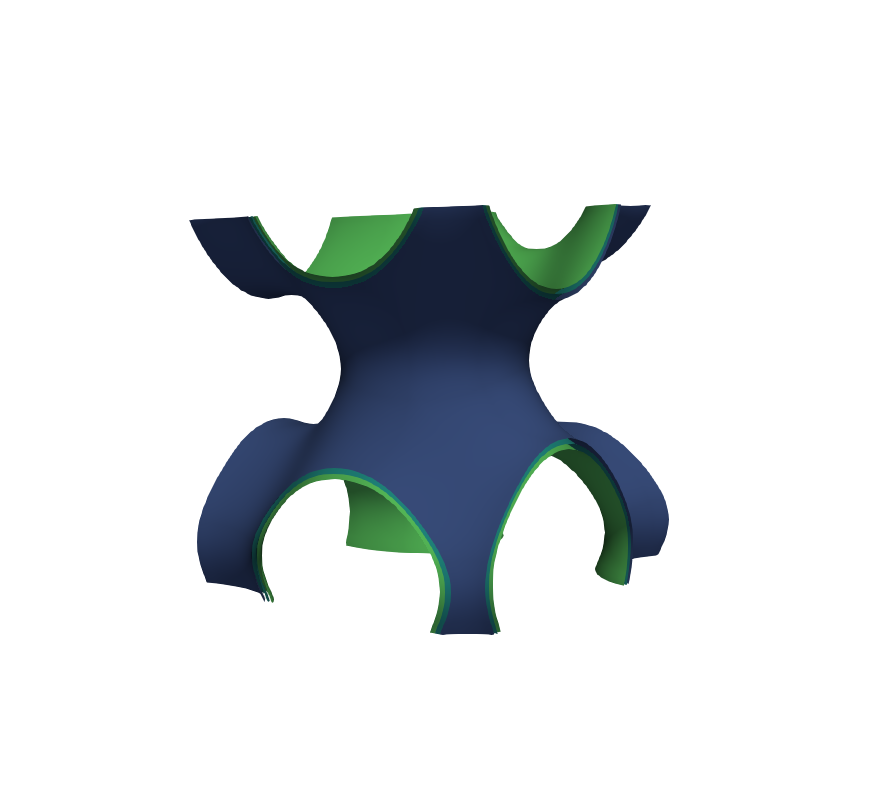}
		\subcaption{steady state(contour)}
		\label{3d_eg1_contour}
	\end{minipage}
	\begin{minipage}[c]{0.30\textwidth}
		\centering
		\includegraphics[width=\textwidth]{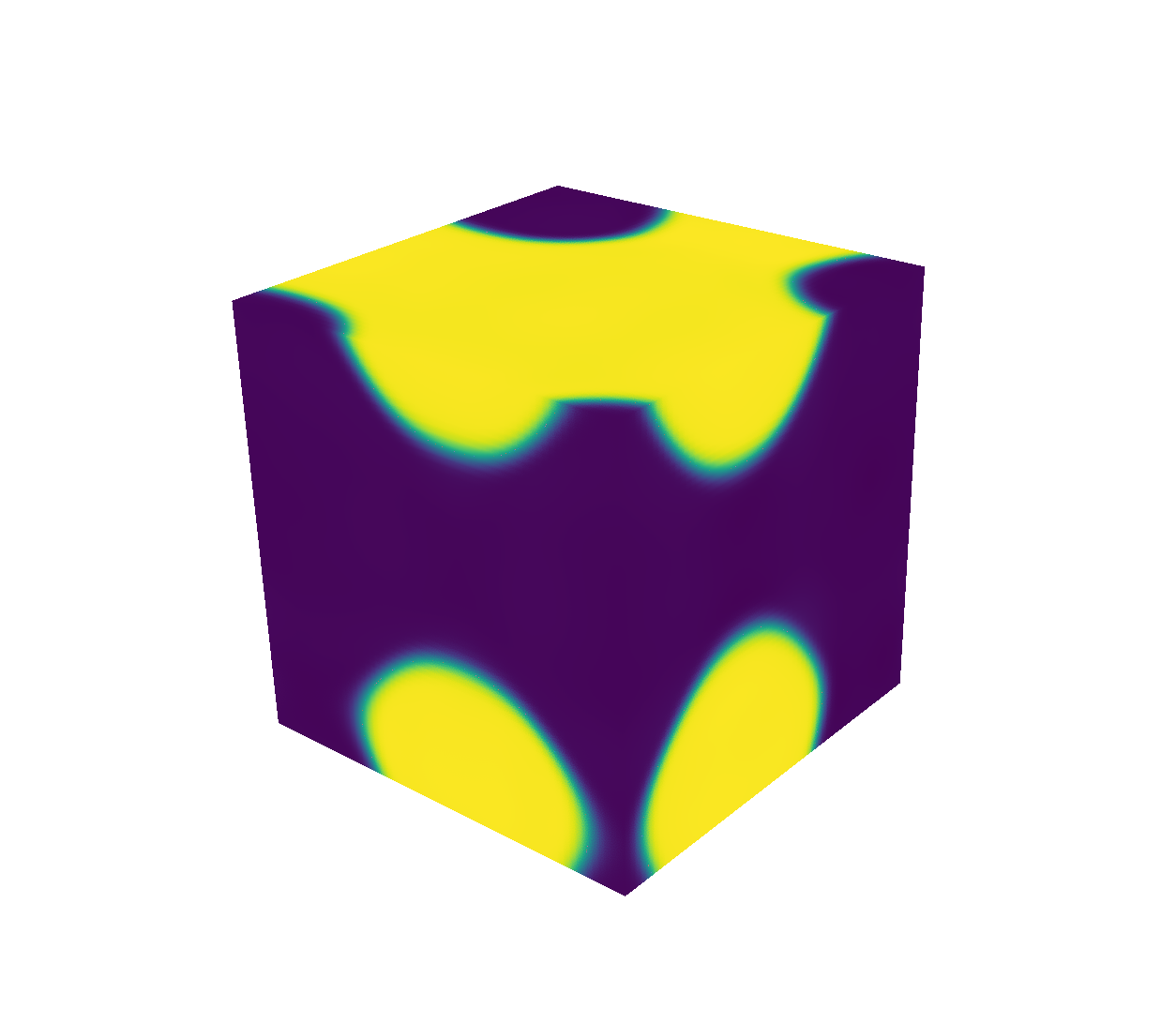}
		\subcaption{steady state(volume)}
		\label{3d_eg1_volume}
	\end{minipage}
    \caption{Case 1 (3D): (a) Random initialization generated by Kaiming initialization. 
(b) The interfacial layer of the droplet-type steady-state solution of the three-dimensional CH equation. 
(c) Volume rendering illustrating the three-dimensional morphology of the phase-separated structure.}
	\label{3d_eg1}
\end{figure}

Figure~\ref{3d_eg1}(a) shows the random initialization using Kaiming initialization, while 

Figure~\ref{3d_eg1}(b) presents the interfacial layer of the corresponding steady-state solution. The associated volume rendering steady-state is shown in Figure~\ref{3d_eg1}(c), which reveals the three-dimensional morphology of the interface. The solution exhibits a clear phase separation structure, $\phi=1$ (yellow part) and $\phi=-1$ (dark part). 
This configuration can be viewed as a three-dimensional extension of the droplet-type structure observed in the two-dimensional case (Figure~\ref{2d_eg5}(b)), as confirmed by its cross-sections at $x=0$ or $y=0$.

\begin{figure}[H]
	\centering
	\begin{minipage}[c]{0.30\textwidth}
		\centering
		\includegraphics[width=\textwidth]{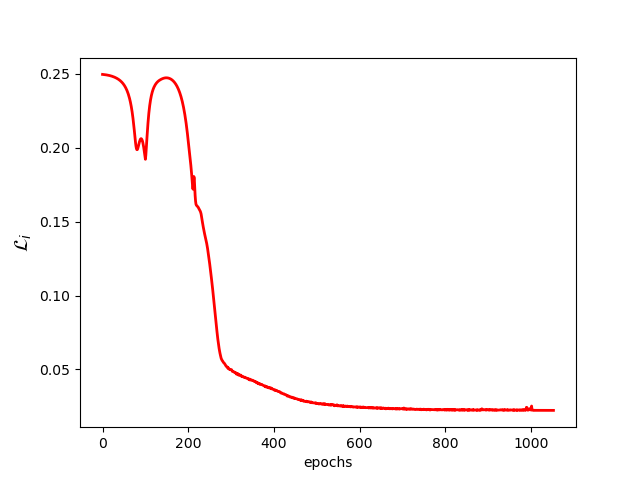}
		\subcaption{energy loss}
		\label{3d_eg1_lossi}
	\end{minipage} 
	\begin{minipage}[c]{0.30\textwidth}
		\centering
		\includegraphics[width=\textwidth]{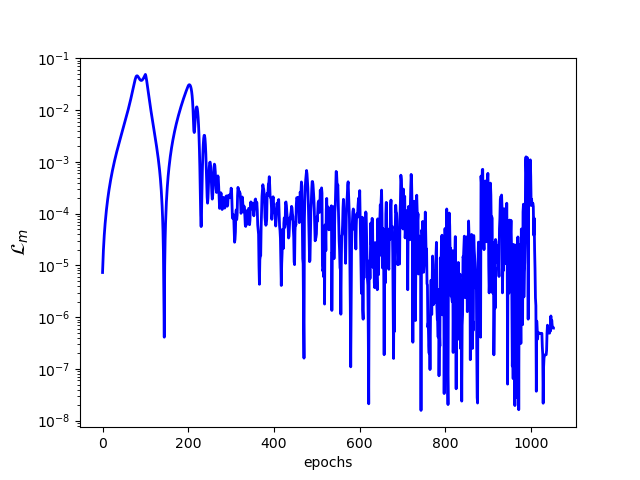}
		\subcaption{mass constraint loss}
		\label{3d_eg1_lossm}
	\end{minipage}
	\begin{minipage}[c]{0.30\textwidth}
		\centering
		\includegraphics[width=\textwidth]{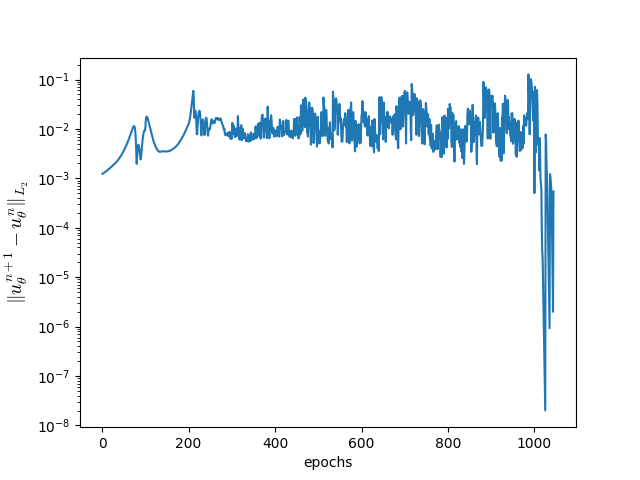}
		\subcaption{ error decay}
		\label{3d_eg1_force}
	\end{minipage}
    \caption{Case 1 (3D): Evolution of (a) energy loss, (b) mass constraint loss, and (c) error during training for the droplet-type steady-state solution.}
	\label{3d_eg1_losses}
\end{figure}

Figures~\ref{3d_eg1_losses}(a) and (b) show the evolution of the energy loss and mass constraint loss, respectively. The steady-state energy is approximately $0.0219$. A mild increase in the energy loss around epoch $1000$ is observed, which is attributed to the transition from the Adam optimizer to L-BFGS for fine-tuning. The error curve in Figure~\ref{3d_eg1_losses}(c) confirms stable convergence.
 
\paragraph{\textit{Case 2}}

\begin{figure}[H]
	\centering
	\begin{minipage}[c]{0.30\textwidth}
		\centering
		\includegraphics[width=\textwidth]{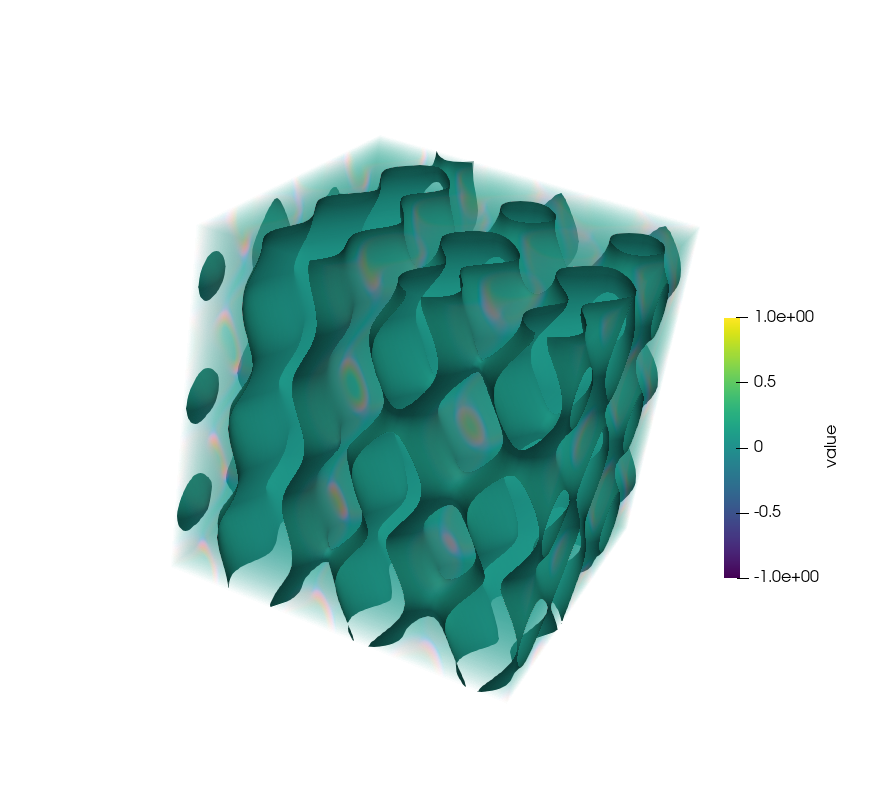}
		\subcaption{random initialization}
		\label{3d_eg2_initial}
	\end{minipage} 
	\begin{minipage}[c]{0.30\textwidth}
		\centering
		\includegraphics[width=\textwidth]{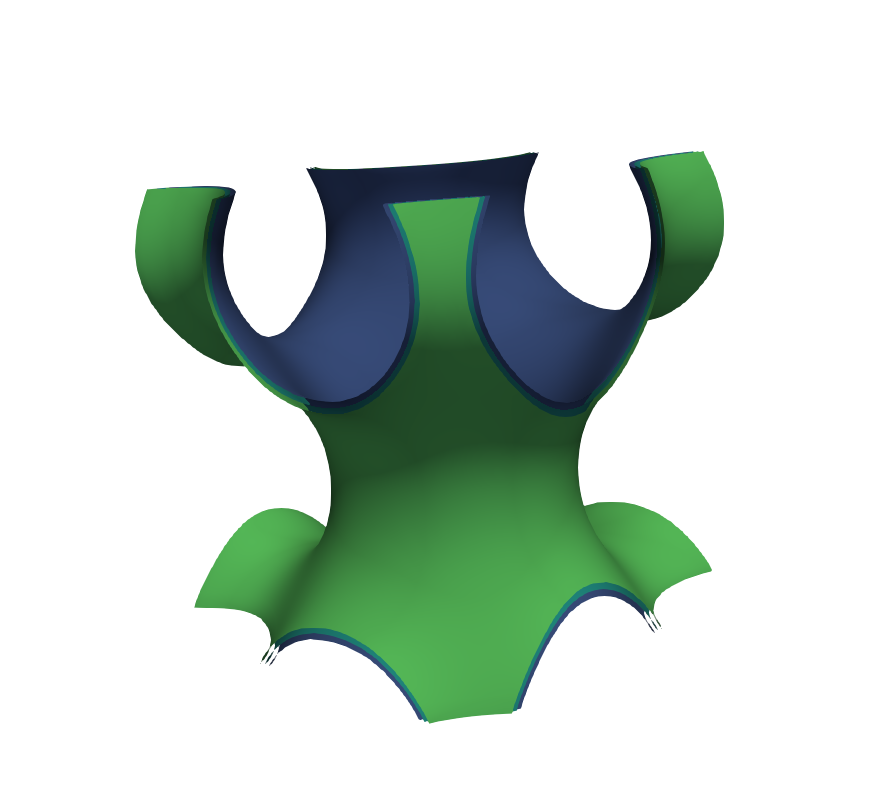}
		\subcaption{steady state(contour)}
		\label{3d_eg2_contour}
	\end{minipage}
	\begin{minipage}[c]{0.30\textwidth}
		\centering
		\includegraphics[width=\textwidth]{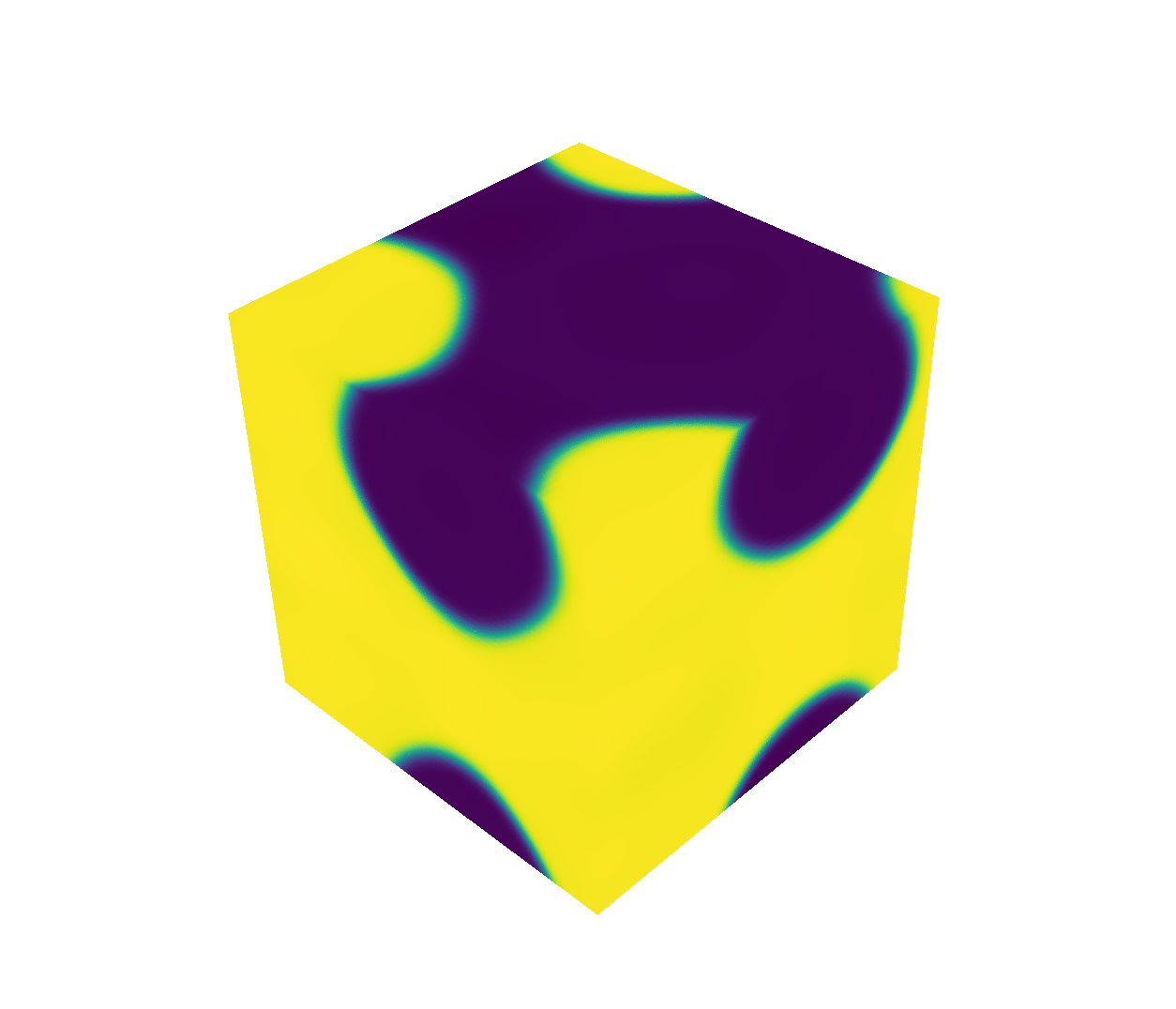}
		\subcaption{steady state(volume)}
		\label{3d_eg2_volume}
	\end{minipage}
    \caption{Case 2 (3D): (a) Random initialization. 
(b) Computed inverted droplet-type steady-state solution, where the two phases are exchanged compared to Case~1. 
(c) Volume rendering showing a geometry similar to Case~1 with reversed phase values.}
	\label{3d_eg2}
\end{figure}

Figure~\ref{3d_eg2}(a) shows another random initialization, and Figure~\ref{3d_eg2}(b) displays the interfacial layer of the resulting steady state. Similar to Case~1, the solution forms a droplet-like phase-separated structure; however, the roles of the two phases are reversed, resulting in an inverted configuration.

The corresponding volume rendering (Figure~\ref{3d_eg2}(c)) shows a geometry nearly identical to that of Case~1, except for the sign inversion of the phases. Cross-sectional views on $x=0$, $y=0$ or $z=0$ again match the corresponding two-dimensional pattern.

\begin{figure}[H]
	\centering
	\begin{minipage}[c]{0.30\textwidth}
		\centering
		\includegraphics[width=\textwidth]{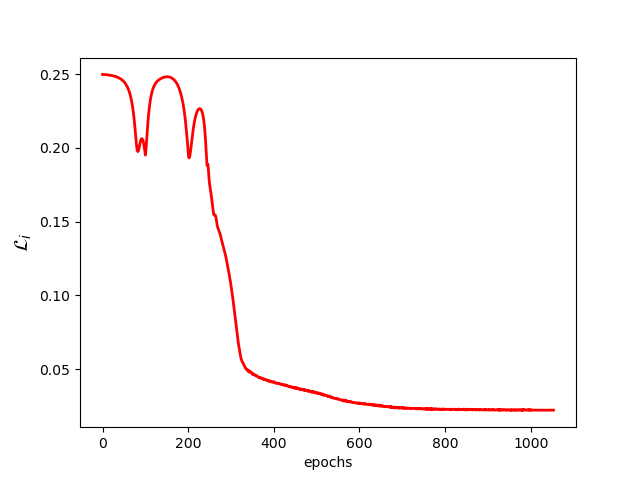}
		\subcaption{energy loss}
		\label{3d_eg2_lossi}
	\end{minipage} 
	\begin{minipage}[c]{0.30\textwidth}
		\centering
		\includegraphics[width=\textwidth]{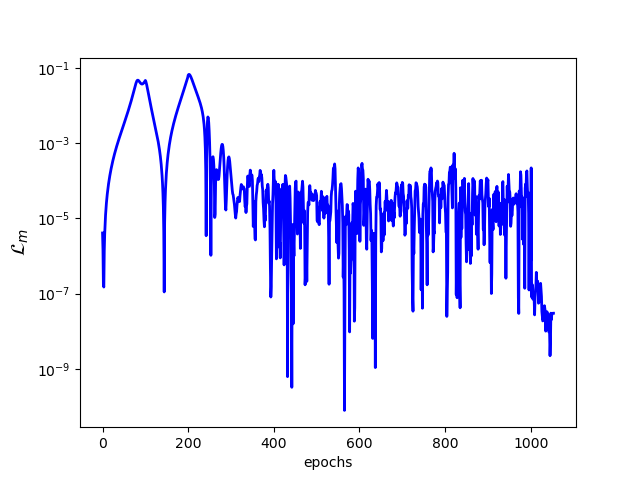}
		\subcaption{mass constraint loss}
		\label{3d_eg2_lossm}
	\end{minipage}
	\begin{minipage}[c]{0.30\textwidth}
		\centering
		\includegraphics[width=\textwidth]{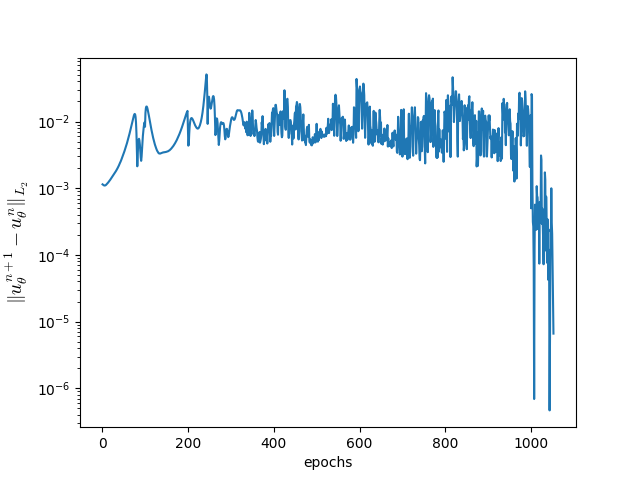}
		\subcaption{ error decay }
		\label{3d_eg2_force}
	\end{minipage}
    \caption{Case 2 (3D): Evolution of (a) energy loss, (b) mass constraint loss, and (c) error during training for the inverted droplet-type solution.}
	\label{3d_eg2_losses}
\end{figure}

The steady-state energy is approximately $0.0219$, indicating that the two configurations correspond to energetically comparable local minimizers. The decay of the energy loss and mass constraint loss is shown in Figures~\ref{3d_eg2_losses}(a) and (b), while the error evolution in Figure~\ref{3d_eg2_losses}(c) further confirms convergence.

\paragraph{\textit{Case 3}}

\begin{figure}[H]
	\centering
	\begin{minipage}[c]{0.45\textwidth}
		\centering
		\includegraphics[width=\textwidth]{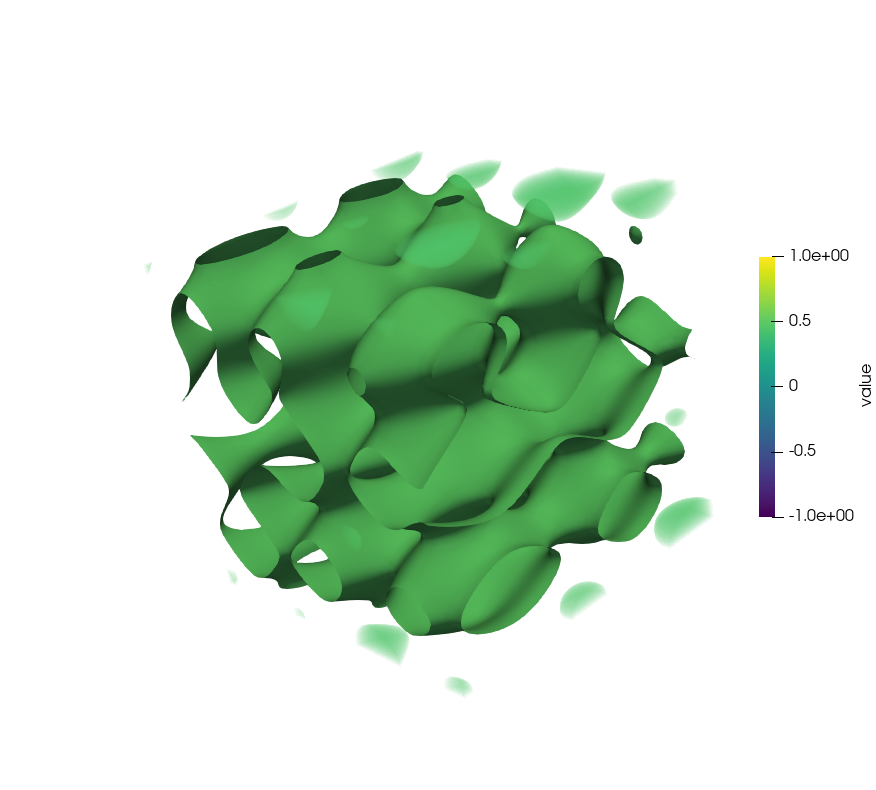}
		\subcaption{random initializion}
		\label{3d_eg3_initial}
	\end{minipage} 
	\begin{minipage}[c]{0.45\textwidth}
		\centering
		\includegraphics[width=\textwidth]{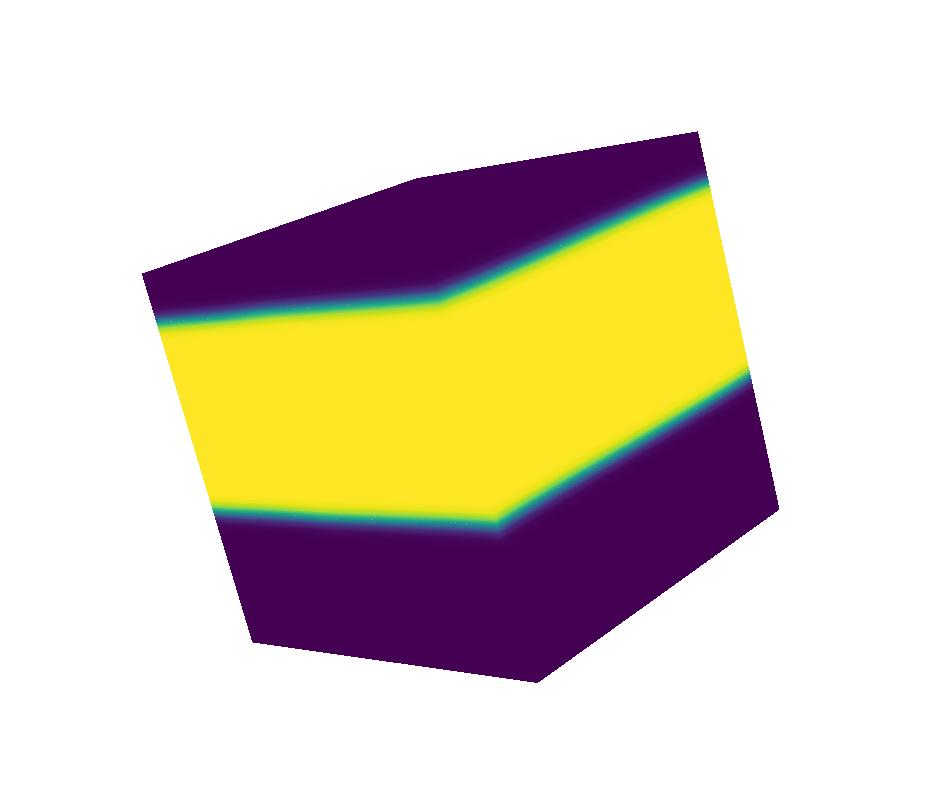}
		\subcaption{steady state(volume)}
		\label{3d_eg3_volume}
	\end{minipage}
    \caption{Case 3 (3D): (a) Random initialization. 
(b) Lamellar ( layered) steady-state solution, representing a three-dimensional extension of the stripe phase in two dimensions.}
	\label{3d_eg3}
\end{figure}

Figures~\ref{3d_eg3}(a) and (b) show the initialization and the computed steady state, respectively. In contrast to the droplet-type structures in Cases~1 and~2, this solution can be interpreted as a three-dimensional extension of the lamellar phase observed in the two-dimensional case. 
The solution exhibits a layered phase separation structure, where the two phases organize into alternating regions across the domain, which is characterized by coherent layered interfaces.
This interpretation is further supported by the cross-sectional slices (e.g., at $x=0$ or $y=0$), which reveal clear stripe-like patterns consistent with the lamellar solutions observed in two dimensions (Figure~\ref{fig_2d_strip}).

\begin{figure}[H]
	\centering
	\begin{minipage}[c]{0.30\textwidth}
		\centering
		\includegraphics[width=\textwidth]{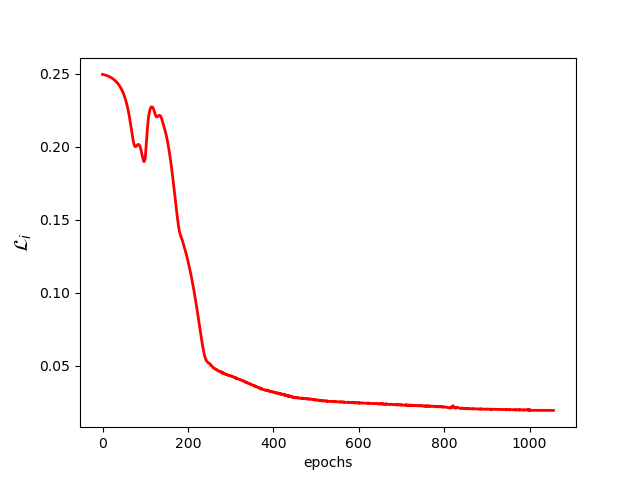}
		\subcaption{energy loss}
		\label{3d_eg3_lossi}
	\end{minipage} 
	\begin{minipage}[c]{0.30\textwidth}
		\centering
		\includegraphics[width=\textwidth]{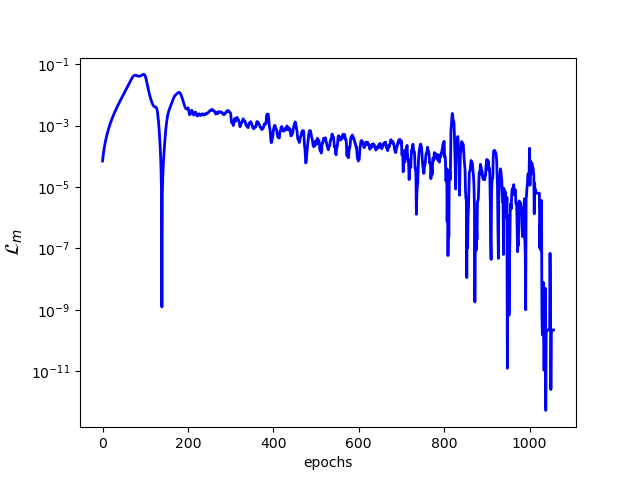}
		\subcaption{mass constraint loss}
		\label{3d_eg3_lossm}
	\end{minipage}
	\begin{minipage}[c]{0.30\textwidth}
		\centering
		\includegraphics[width=\textwidth]{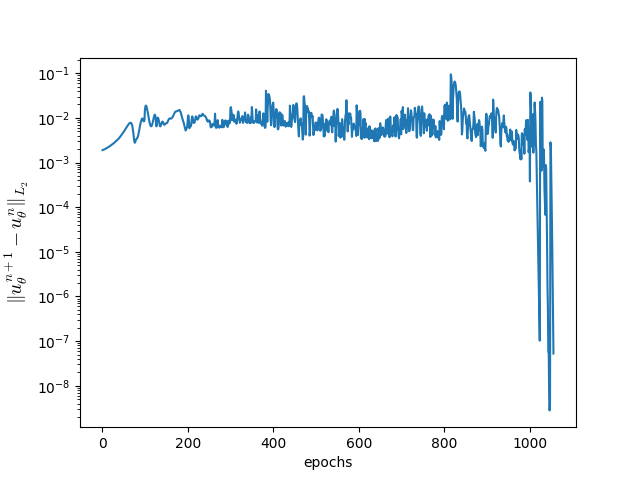}
		\subcaption{ error decay }
		\label{3d_eg3_force}
	\end{minipage}
	\caption{Case 3 (3D): Evolution of (a) energy loss, (b) mass constraint loss, and (c) error during training for the lamellar steady-state solution.}
	\label{3d_eg3_losses}
\end{figure}

The steady-state energy is approximately $0.0187$, which is lower than that of the droplet-type solutions, suggesting a more energetically favorable configuration. Figures~\ref{3d_eg3_losses}(a) and (b) show the decay of the energy loss and mass constraint loss, respectively, while the error curve in Figure~\ref{3d_eg3_losses}(c) demonstrates stable convergence.

\paragraph{\textit{Case 4}}

\begin{figure}[H]
	\centering
	\begin{minipage}[c]{0.30\textwidth}
		\centering
		\includegraphics[width=\textwidth]{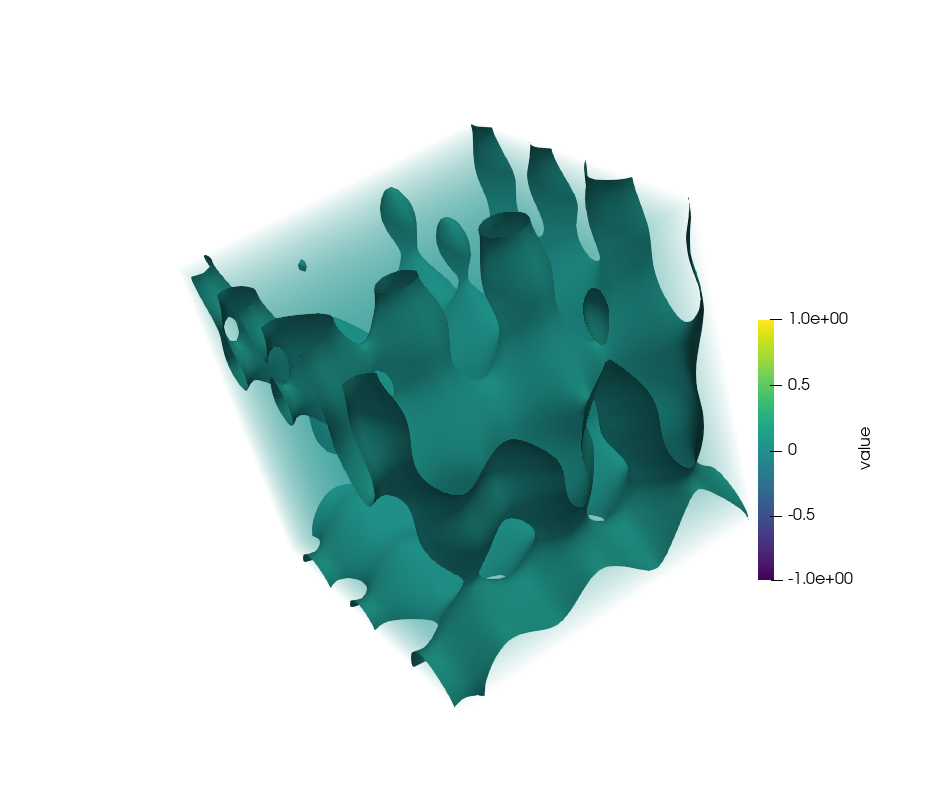}
		\subcaption{}
		%\label{3d_eg4_initial}
	\end{minipage}%
	\begin{minipage}[c]{0.30\textwidth}
		\includegraphics[width=\textwidth]{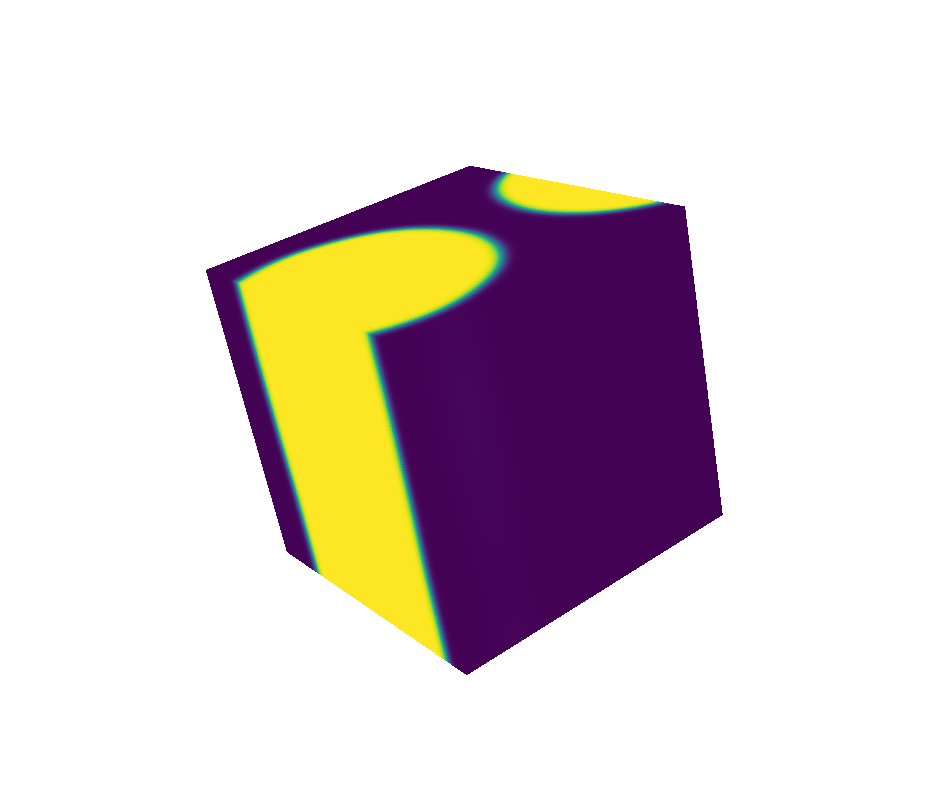}
		\subcaption{}
		%\label{3d_eg4_volume}
	\end{minipage}%    
	\begin{minipage}[c]{0.30\textwidth}
		\includegraphics[width=\textwidth]{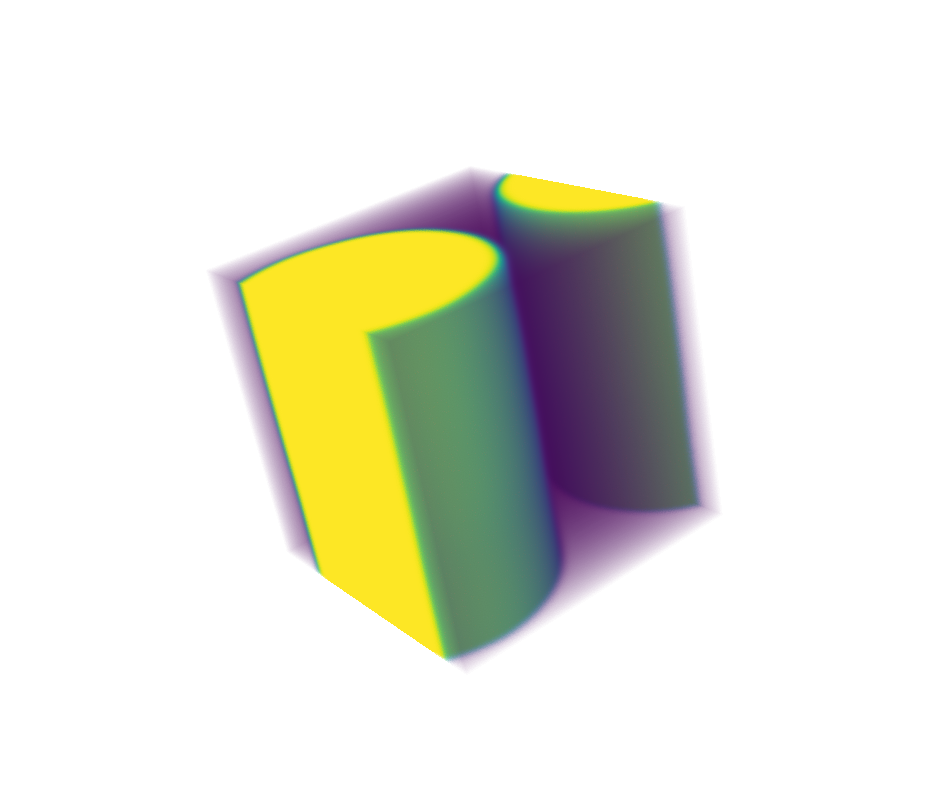}
		\subcaption{}
		%\label{3d_eg4_initial}
	\end{minipage}% 
	% \begin{minipage}[c]{0.25\textwidth}
	% 	\includegraphics[width=\textwidth]{pictures/3dim/5/3dpred4.png}
	% 	\subcaption{}
	% 	%\label{3d_eg4_volume}
	% \end{minipage}    
    \caption{Case 4 (3D): (a) Random initialization. 
(b)--(c) the same tubular (cylindrical) steady-state solution with different opacity settings for better observation of the structure.}
	\label{3d_eg5}
\end{figure}

Figures~\ref{3d_eg5}(a) and (b) illustrate the random initialization and the corresponding steady-state solution, respectively. The same steady state is further visualized in Figure~\ref{3d_eg5}(c) and (d) using different opacity settings, which provide a clearer view of the underlying three-dimensional structure.
The solution exhibits a cylindrical (or tubular) morphology, where one phase ($\phi\approx 1$) is embedded within the surrounding phase ($\phi\approx -1$).
Cross-sectional profiles reveal distinct lower-dimensional characteristics: the slice at $z=0$ corresponds to a droplet-type structure in two dimensions (see Figure~\ref{fig_droplet_2}), while slices taken at $x=x_0$ or $y=y_0$ display lamellar (striped) patterns or nearly homogeneous states, depending on the location. These observations highlight the anisotropic and spatially heterogeneous nature of the solution structure.

\begin{figure}[H]
	\centering
	\begin{minipage}[c]{0.30\textwidth}
		\centering
		\includegraphics[width=\textwidth]{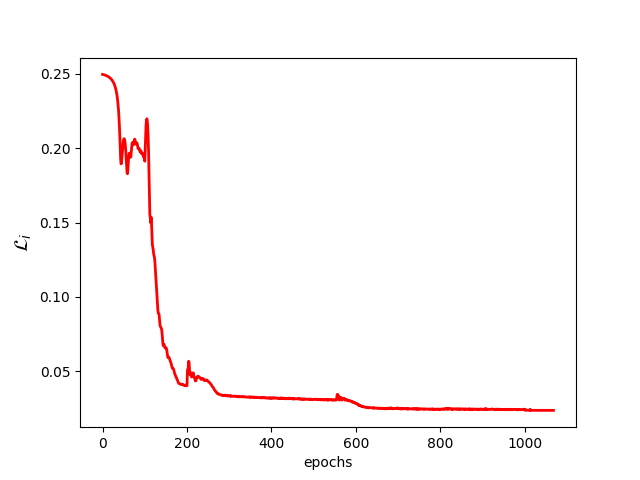}
		\subcaption{energy loss}
		\label{3d_eg4_lossi}
	\end{minipage} 
	\begin{minipage}[c]{0.30\textwidth}
		\centering
		\includegraphics[width=\textwidth]{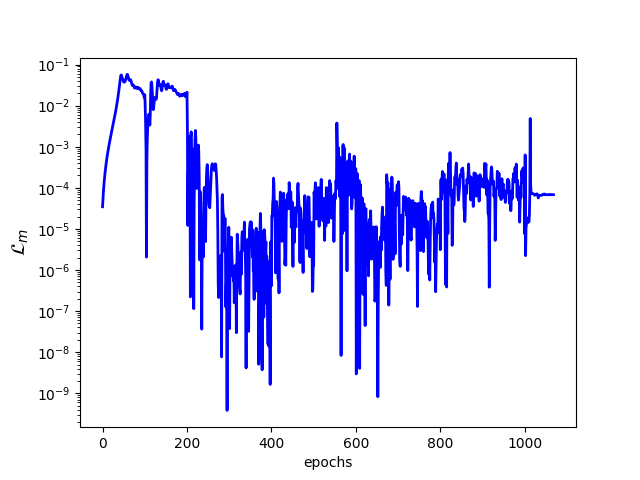}
		\subcaption{mass constraint loss}
		\label{3d_eg4_lossm}
	\end{minipage}
	\begin{minipage}[c]{0.30\textwidth}
		\centering
		\includegraphics[width=\textwidth]{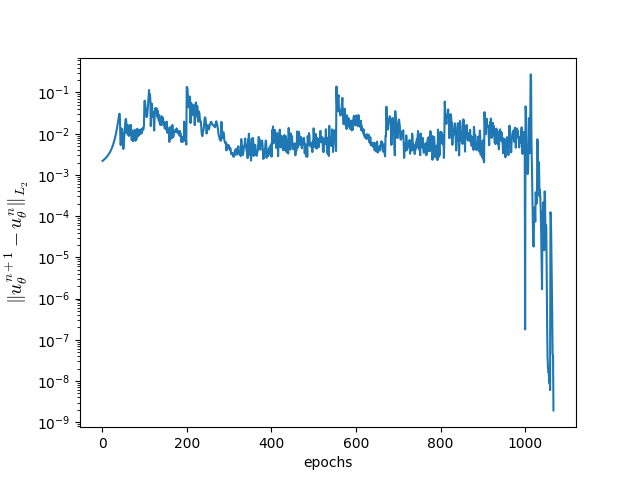}
		\subcaption{error decay}
		\label{3d_eg4_force}
	\end{minipage}
    \caption{Case 4 (3D): Evolution of (a) energy loss, (b) mass constraint loss, and (c) error during training for the tubular steady-state solution.}
	\label{3d_eg5_losses}
\end{figure}

The steady-state energy is approximately $0.0238$. Figures~\ref{3d_eg5_losses}(a) and (b) show the evolution of the energy loss and mass constraint loss, respectively, while Figure~\ref{3d_eg5_losses}(c) presents the error decay, which confirms the convergence of the algorithm. %As in Case~1, a mild fluctuation in the energy loss is observed near epoch $1000$, which is attributed to the transition from the Adam optimizer to L-BFGS.
The computational time required to reach convergence for these three-dimensional cases is approximately 95 seconds, demonstrating the efficiency of the proposed method.

These results demonstrate that the proposed framework is capable of capturing a rich variety of three-dimensional phase separation patterns, including droplet-type, lamellar, and tubular structures, highlighting its effectiveness and its strong ability in exploring complex high-dimensional energy landscapes.

% \begin{figure}[H]
% 	\centering
% 	%\begin{minipage}[c]{0.45\textwidth}
% 		\centering
% 		\includegraphics[width=0.8\textwidth]{pictures/3dim/5/initial.png}
%         \includegraphics[width=0.8\textwidth]{pictures/3dim/5/3dpred2.png}
%         \includegraphics[width=0.8\textwidth]{pictures/3dim/5/3dpred3.png}
% 		\subcaption{random initialization}
% 		%\label{3d_eg4_initial}
% 	%\end{minipage} 
%     \caption{3D, Case 4: (a) Random initialization. 
% (b) Tubular (cylindrical) steady-state solution.}%, where one phase is embedded in the other phase. }
% 	\label{3d_eg5}
% \end{figure}

\section{Conclusion}\label{Conclusion}

In this work, we propose a deep learning-based framework, namely the Deep Ritz method, for computing steady-state solutions of the Cahn--Hilliard equation under periodic boundary conditions. The method effectively alleviates the curse of dimensionality encountered in traditional numerical approaches, particularly in high-dimensional settings.
To enforce the mass conservation constraint, we incorporate an enhanced augmented Lagrangian formulation, which ensures both accuracy and stability during training. In addition, the random/separable Fourier feature mapping is employed to naturally satisfy periodic boundary conditions, while significantly improving the expressive capability of the neural network for capturing oscillatory and nontrivial solution structures.

Extensive numerical experiments in one-, two-, and three-dimensional cases demonstrate that the proposed method is capable of efficiently identifying multiple steady states, including both trivial and non-trivial solutions. In particular, a rich variety of phase-separated structures—such as droplet-type, lamellar, and tubular patterns—are successfully captured. These results highlight the robustness and effectiveness of the proposed framework in exploring complex high-dimensional energy landscapes and locating diverse local minimizers of the underlying energy functional.
Compared with conventional grid-based methods, the proposed approach offers notable advantages in terms of flexibility, scalability, and computational efficiency, especially for high-dimensional problems. Moreover, the method does not rely on carefully designed initial guesses and is able to discover multiple solution branches through different random initializations.

The proposed framework can be further extended to a broader class of phase-field models for steady-state computation. In addition, it provides a promising foundation for locating transition states by combining with saddle point search algorithms, which will be investigated in our future work.

\section*{Acknowledgments}
 Shuting GU acknowledges the support of NSFC 12571465 and Guangdong Basic and Applied Basic Research Foundation 2026A1515012326.
 
%%Vancouver style references.
\bibliographystyle{model1-num-names}
\bibliography{refs}

@article{cahn1958free,
  title={Free energy of a nonuniform system. I. Interfacial free energy},
  author={Cahn, John W and Hilliard, John E},
  journal={The Journal of Chemical Physics},
  volume={28},
  number={2},
  pages={258--267},
  year={1958},
  publisher={AIP}
}

@book{chaikin2000principles,
  title={Principles of condensed matter physics},
  author={Chaikin, Paul M and Lubensky, Thomas C},
  year={2000},
  publisher={Cambridge University Press}
}

@article{elliott1989global,
  title={On the Cahn–Hilliard equation},
  author={Elliott, Charles M and Zheng, Songmu},
  journal={Archive for Rational Mechanics and Analysis},
  volume={96},
  number={4},
  pages={339--357},
  year={1989},
  publisher={Springer}
}

@article{raissi2019physics,
  title={Physics-informed neural networks: A deep learning framework for solving forward and inverse problems involving nonlinear partial differential equations},
  author={Raissi, Maziar and Perdikaris, Paris and Karniadakis, George E},
  journal={Journal of Computational Physics},
  volume={378},
  pages={686--707},
  year={2019},
  publisher={Elsevier}
}

@article{weinan2018deep,
  title={The deep Ritz method: A deep learning-based numerical algorithm for solving variational problems},
  author={E, Weinan and Yu, Bing},
  journal={Communications in Mathematics and Statistics},
  volume={6},
  number={1},
  pages={1--12},
  year={2018},
  publisher={Springer}
}

@inproceedings{Lu2021Generalization,
  title = {A Priori Generalization Analysis of the Deep Ritz Method for Solving High Dimensional Elliptic Partial Differential Equations},
  author = {Yulong Lu and Jianfeng Lu and Min Wang},
  booktitle = {Proceedings of Thirty Fourth Conference on Learning Theory},
  pages = {3196--3241},
  year = {2021},
  volume = {134},
  series = {Proceedings of Machine Learning Research}
}

@InProceedings{He_2015_ICCV,
  author    = {He, Kaiming and Zhang, Xiangyu and Ren, Shaoqing and Sun, Jian},
  title     = {Delving Deep into Rectifiers: Surpassing Human-Level Performance on ImageNet Classification},
  booktitle = {Proceedings of the IEEE International Conference on Computer Vision (ICCV)},
  year      = {2015},
  pages     = {1026--1034},
  month     = {Dec},
  doi       = {10.1109/ICCV.2015.123},
  url       = {https://openaccess.thecvf.com/content_iccv_2015/html/He_Delving_Deep_into_ICCV_2015_paper.html}
}

@inproceedings{Tancik2020FourierFeatures,
  title     = {Fourier Features Let Networks Learn High Frequency Functions in Low Dimensional Domains},
  author    = {Tancik, Matthew and Srinivasan, Pratul P. and Mildenhall, Ben and Fridovich-Keil, Sara and Raghavan, Nithin and Singhal, Utkarsh and Ramamoorthi, Ravi and Barron, Jonathan T. and Ng, Ren},
  booktitle = {Advances in Neural Information Processing Systems (NeurIPS)},
  year      = {2020},
  url       = {https://arxiv.org/abs/2006.10739}
}

@article{Rahimi2007RandomFeatures,
  title   = {Random Features for Large-Scale Kernel Machines},
  author  = {Rahimi, Ali and Recht, Benjamin},
  journal = {Advances in Neural Information Processing Systems (NeurIPS)},
  year    = {2007},
  pages   = {1177--1184}
}

@article{han2018solving,
  title={Solving high-dimensional partial differential equations using deep learning},
  author={Han, Jiequn and Jentzen, Arnulf and E, Weinan},
  journal={Proceedings of the National Academy of Sciences},
  volume={115},
  number={34},
  pages={8505--8510},
  year={2018},
  publisher={National Academy of Sciences}
}

@article{han2025brief,
  title={A brief review of the Deep BSDE method for solving high-dimensional partial differential equations},
  author={Han, Jiequn and Jentzen, Arnulf and others},
  journal={arXiv preprint arXiv:2505.17032},
  year={2025}
}

@article{raissi2017physics,
  title={Physics informed deep learning (part i): Data-driven solutions of nonlinear partial differential equations},
  author={Raissi, Maziar and Perdikaris, Paris and Karniadakis, George Em},
  journal={arXiv preprint arXiv:1711.10561},
  year={2017}
}

@article{sirignano2018dgm,
  title={DGM: A deep learning algorithm for solving partial differential equations},
  author={Sirignano, Justin and Spiliopoulos, Konstantinos},
  journal={Journal of computational physics},
  volume={375},
  pages={1339--1364},
  year={2018},
  publisher={Elsevier}
}

@incollection{raissi2024forward,
  title={Forward--backward stochastic neural networks: deep learning of high-dimensional partial differential equations},
  author={Raissi, Maziar},
  booktitle={Peter Carr Gedenkschrift: Research Advances in Mathematical Finance},
  pages={637--655},
  year={2024},
  publisher={World Scientific}
}

@article{zhang2022fbsde,
  title={FBSDE based neural network algorithms for high-dimensional quasilinear parabolic PDEs},
  author={Zhang, Wenzhong and Cai, Wei},
  journal={Journal of Computational Physics},
  volume={470},
  pages={111557},
  year={2022},
  publisher={Elsevier}
}

@article{zang2020weak,
  title={Weak adversarial networks for high-dimensional partial differential equations},
  author={Zang, Yaohua and Bao, Gang and Ye, Xiaojing and Zhou, Haomin},
  journal={Journal of Computational Physics},
  volume={411},
  pages={109409},
  year={2020},
  publisher={Elsevier}
}

@InProceedings{pmlr-v80-long18a,
  title = 	 {{PDE}-Net: Learning {PDE}s from Data},
  author =       {Long, Zichao and Lu, Yiping and Ma, Xianzhong and Dong, Bin},
  booktitle = 	 {Proceedings of the 35th International Conference on Machine Learning},
  pages = 	 {3208--3216},
  year = 	 {2018},
  editor = 	 {Dy, Jennifer and Krause, Andreas},
  volume = 	 {80},
  series = 	 {Proceedings of Machine Learning Research},
  month = 	 {10--15 Jul},
  publisher =    {PMLR},
  url = 	 {https://proceedings.mlr.press/v80/long18a.html}
}

@article{long2019pde,
  title={PDE-Net 2.0: Learning PDEs from data with a numeric-symbolic hybrid deep network},
  author={Long, Zichao and Lu, Yiping and Dong, Bin},
  journal={Journal of Computational Physics},
  volume={399},
  pages={108925},
  year={2019},
  publisher={Elsevier}
}

@article{hure2020deep,
  title={Deep backward schemes for high-dimensional nonlinear PDEs},
  author={Hur{\'e}, C{\^o}me and Pham, Huy{\^e}n and Warin, Xavier},
  journal={Mathematics of Computation},
  volume={89},
  number={324},
  pages={1547--1579},
  year={2020}
}

@article{germain2022approximation,
  title={Approximation error analysis of some deep backward schemes for nonlinear PDEs},
  author={Germain, Maximilien and Pham, Huyen and Warin, Xavier},
  journal={SIAM Journal on Scientific Computing},
  volume={44},
  number={1},
  pages={A28--A56},
  year={2022},
  publisher={SIAM}
}

@article{beck2021deep,
  title={Deep splitting method for parabolic PDEs},
  author={Beck, Christian and Becker, Sebastian and Cheridito, Patrick and Jentzen, Arnulf and Neufeld, Ariel},
  journal={SIAM Journal on Scientific Computing},
  volume={43},
  number={5},
  pages={A3135--A3154},
  year={2021},
  publisher={SIAM}
}

@article{ji2020three,
  title={Three algorithms for solving high-dimensional fully coupled FBSDEs through deep learning},
  author={Ji, Shaolin and Peng, Shige and Peng, Ying and Zhang, Xichuan},
  journal={IEEE Intelligent Systems},
  volume={35},
  number={3},
  pages={71--84},
  year={2020},
  publisher={IEEE}
}

@article{han2016deep,
  title={Deep learning approximation for stochastic control problems},
  author={Han, Jiequn and others},
  journal={arXiv preprint arXiv:1611.07422},
  year={2016}
}

@article{cai2023deepmartnet,
  title={DeepMartNet--A Martingale based Deep Neural Network Learning Algorithm for Eigenvalue/BVP Problems and Optimal Stochastic Controls},
  author={Cai, Wei},
  journal={arXiv preprint arXiv:2307.11942},
  year={2023}
}

@article{han2017deep,
  title={Deep learning-based numerical methods for high-dimensional parabolic partial differential equations and backward stochastic differential equations},
  author={Han, Jiequn and Jentzen, Arnulf and others},
  journal={Communications in mathematics and statistics},
  volume={5},
  number={4},
  pages={349--380},
  year={2017},
  publisher={Springer}
}

@article{beck2019machine,
  title={Machine learning approximation algorithms for high-dimensional fully nonlinear partial differential equations and second-order backward stochastic differential equations},
  author={Beck, Christian and E, Weinan and Jentzen, Arnulf},
  journal={Journal of Nonlinear Science},
  volume={29},
  number={4},
  pages={1563--1619},
  year={2019},
  publisher={Springer}
}

@article{zhang2022deep,
  title={Deep network approximation: Achieving arbitrary accuracy with fixed number of neurons},
  author={Zhang, Shijun and Shen, Zuowei and Yang, Haizhao},
  journal={Journal of Machine Learning Research},
  volume={23},
  number={276},
  pages={1--60},
  year={2022}
}

@inproceedings{rahaman2019spectral,
  title={On the spectral bias of neural networks},
  author={Rahaman, Nasim and Baratin, Aristide and Arpit, Devansh and Draxler, Felix and Lin, Min and Hamprecht, Fred and Bengio, Yoshua and Courville, Aaron},
  booktitle={International conference on machine learning},
  pages={5301--5310},
  year={2019},
  organization={PMLR}
}

@article{tancik2020fourier,
  title={Fourier features let networks learn high frequency functions in low dimensional domains},
  author={Tancik, Matthew and Srinivasan, Pratul and Mildenhall, Ben and Fridovich-Keil, Sara and Raghavan, Nithin and Singhal, Utkarsh and Ramamoorthi, Ravi and Barron, Jonathan and Ng, Ren},
  journal={Advances in neural information processing systems},
  volume={33},
  pages={7537--7547},
  year={2020}
}

@article{hornik1991approximation,
  title={Approximation capabilities of multilayer feedforward networks},
  author={Hornik, Kurt},
  journal={Neural networks},
  volume={4},
  number={2},
  pages={251--257},
  year={1991},
  publisher={Elsevier}
}

\end{document}